\documentclass[final,1p]{elsarticle}

\usepackage{longtable}

\journal{Computer-Aided Design}

\usepackage{bm}

\usepackage{tikz}

\usepackage{amsmath,amssymb}
\usepackage{amsfonts}
\usepackage{pmbb-sym}
\usepackage{here}

\biboptions{numbers,sort&compress} 

\def\diff{\mathrm{d}}

\usepackage{subfigure}
\def\pp#1#2{\frac{\partial #1}{\partial #2}}

\def\abs#1{\left|#1\right|}

\def\Re{\mathop{\mathrm{Re}}}

\def\imath{\mathrm{i}}

\usepackage{siunitx}
\NewDocumentCommand\numprint{m}{\num[round-mode = places]{#1}}
\NewDocumentCommand\nprounddigits{m}{\sisetup{round-precision = #1}}
\def\npproductsign#1{} 

\usepackage{braket}

\usepackage[colorlinks,bookmarksopen,bookmarksnumbered,citecolor=red,urlcolor=red,breaklinks=true]{hyperref} 
\def\equationautorefname~#1\null{%
  Eq.~(#1)\null
}
\def\algorithmautorefname~#1\null{%
  Algorithm~#1\null
}
\def\sectionautorefname~#1\null{%
  Section~#1\null
}
\def\subsectionautorefname~#1\null{%
  Subsection~#1\null
}
\def\subsubsectionautorefname~#1\null{%
  Subsubsection~#1\null
}
\def\tableautorefname~#1\null{%
  Table~#1\null
}
\def\figureautorefname~#1\null{%
  Figure~#1\null
}
\def\subfigureautorefname~#1\null{%
  Figure~#1\null
}

\usepackage{listings}
\usepackage{framed} 

\definecolor{darkgreen}{cmyk}{1,0,1,0.608}
\definecolor{darkorange}{rgb}{1.0, 0.5, 0.2}

\def\ps{{p_s}} 
\def\pt{{p_t}} 
\def\ns{{n_s}}
\def\nt{{n_t}}
\def\NN{N}

\def\norm#1{\left|#1\right|}
\def\uin{u^{\rm in}}

\usepackage{pgfplots}
\usepackage{booktabs}

\usepackage{siunitx}
\graphicspath{{./figure/}}

\iffalse
\renewenvironment{quote}{\begin{framed}%
    \noindent
  }{\end{framed}}
\else

\fi

\begin{document}

\begin{frontmatter}

  \title{A shape optimisation with the isogeometric boundary element method and adjoint variable method for the three-dimensional Helmholtz equation}

  \tnoteref{kaken}
  \tnotetext[kaken]{This work was partially by JSPS KAKENHI Grant Number 18K11335.}
  
  \author[NU]{Toru Takahashi\corref{cor}}
  \ead{toru.takahashi@mae.nagoya-u.ac.jp}
  \cortext[cor]{Corresponding author}
  \author[NU]{Daisuke Sato}
  \author[NU]{Hiroshi Isakari}
  \author[NU]{Toshiro Matsumoto}
  
  \address[NU]{Department of Mechanical Systems Engineering, Graduate School of Engineering, Nagoya University,\\ Furo-cho, Nagoya, Aichi, 464-8603 Japan}

  \begin{abstract}
    This paper presents a shape optimisation system to design the shape of an acoustically-hard object in the three-dimensional open space. Boundary element method (BEM) is suitable to analyse such an exterior field. However, the conventional BEM, which is based on piecewise polynomial shape and interpolation functions, can require many design variables because they are usually chosen as a part of the nodes of the underlying boundary element mesh. In addition, it is not easy for the conventional method to compute the gradient of the sound pressure on the surface, which is necessary to compute the shape derivative of our interest, of a given object. To overcome these issues, we employ the isogeometric boundary element method (IGBEM), which was developed in our previous work. With using the IGBEM, we can design the shape of surfaces through control points of the NURBS surfaces of the target object. We integrate the IGBEM with the nonlinear programming software through the adjoint variable method (AVM), where the resulting adjoint boundary value problem can be also solved by the IGBEM with a slight modification. The numerical verification and demonstration validate our shape optimisation framework.
  \end{abstract}
  
  \begin{keyword}
    Boundary Element Method \sep Isogeometric Analysis \sep Shape Optimisation \sep Adjoint Variable Method \sep Nonlinear Programming Problem
  \end{keyword}
  
\end{frontmatter}

\section{Background and purpose}\label{s:intro}

Interaction of waves with materials in a specific shape/topology can bring exotic wave phenomena such as the extraordinary transmission through a sub-wavelength aperture~\cite{Ebbesen_1998} and the emergence of a collimated beam by corrugating the aperture~\cite{Lezec_2002, Christensen_2007, Zhou_2010, takahashi2014}. In particular, structures with certain periodic patterns, that is, metamaterials have been recently and intensively studied in science and engineering~\cite{liu2011metamaterials,wang2020tunable}.

Shape optimisation is useful in appropriately designing (meta)materials in a wave field of interest. To analyse the wave problem numerically, boundary element method (BEM) is suitable because it can deal with the infinite domain without any absorbing boundary condition, which is necessary for domain-type solvers such as finite element method (FEM) and finite difference method. In addition, boundary-only models that BEM handles fit in shape optimisation, which concerns the deformation of the surface (or boundary) of a target material rather than its inside.

\begin{figure}[H]
  \centering
  \begin{tikzpicture}[scale=.8]
  \begin{axis}[
    xlabel={Year},
    ylabel style={align=center},
    ylabel={Publications},
    xmin=1985, xmax=2025, xtick={1990,2000,...,2020},
    ymin=0, ymax=25, ytick={0,5,...,25},
    xmajorgrids=true,
    ymajorgrids=true,
    legend pos=north west,
    height=150pt,
    width=400pt,
    ybar, bar width=3pt,
    /pgf/number format/.cd, use comma, 1000 sep={}
    ]
    \addplot table [x=Year, y=Publications, ybar] {shape_optimization_boundary_element_method.txt};
    \addlegendentry{``shape optimization'' + ``boundary elemenet method''}
    \addplot table [x=Year, y=Publications] {isogeometric_boundary_element_method.txt};
    \addlegendentry{``isogeomtric boundary elemenet method''}
  \end{axis}
\end{tikzpicture}
  \caption{Publications with the terms both ``shape optimization'' and ``boundary element method'' (coloured in blue) and those with the term ``isogeometric boundary element method'' (in red). The data was obtained from Web of Science on Apr 28, 2021.}
  \label{fig:publications}
\end{figure}
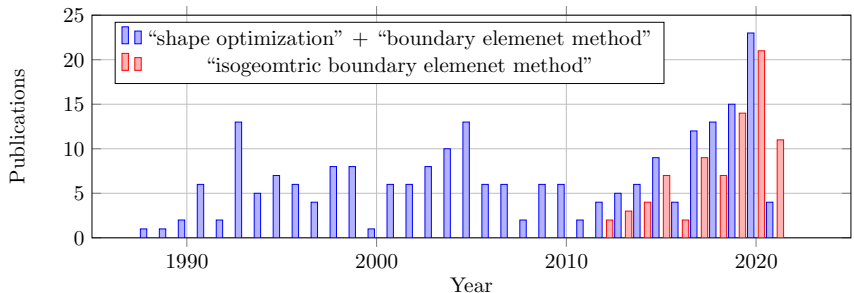

The first study on shape optimisation with using BEM was conducted by Soares et al. in 1984~\cite{soares1984}, which investigated linear-elastostatic problems in 2D. Since then, there are over 200 publications that are related to both shape optimisation and BEM, as shown in Figure~\ref{fig:publications}. The BEM-based shape optimisation is currently being promoted by a new type of BEM, that is, isogeometric BEM (IGBEM), whose number of publications is also shown in the same figure. The IGBEM is characterised by employing the NURBS (including B-spline) function as both shape and interpolation functions, following the concept of isogeometric analysis (IGA)~\cite{hughes2005,cottrell2009}. In this case, one can design the shape of interest through the control points (CPs) associated with the NURBS surface(s). On the other hand, one needs to regard (a part of) the nodes of a boundary element mesh as the design variables in the case of the conventional BEMs, which are based on a piecewise polynomial basis. This can increase the number of design variables unnecessarily, in particular, when the shape of the boundary is complicated and, thus, the mesh is fine. In the IGA, the technique of knot insertion can readily resolve the dilemma between reducing the number of design variables and increasing the resolution of the boundary element analysis. This is the main advantage of the IGBEM over the conventional BEMs, although the formulation and implementation of the former are hard than those of the latter.

Another merit of the IGBEM is that we can easily compute the gradient of the sound pressure at any points (generally except for its boundary, where the other surfaces are connected) on the surface of a scatterer. This is because the sound pressure is usually differentiable over a NURBS surface. This property of the IGBEM is useful when we compute the shape derivative of our interest (see (\ref{eq:sd})). On contrary, the gradient can be discontinuous on the edges of the boundary element mesh in the convectional BEM.

So far, shape optimisations based on the IGBEM have been investigated in terms of potential problems (or steady-state heat problems)~\cite{yoon2015,kostas2015,gillebaart2016,kostas2017,kostas2018}, elastostatic problems~\cite{li2011,lian2016,lian2017,sun2018,li2019,sun2020}, including 2D thermoelastic problem~\cite{yoon2020}, and acoustic problems in concern~\cite{liu2017,takahashi2019ewco,ummidivarapu2020,shaaban2020,shaaban2020b,wang2020,chen2019}. In regard to 2D, Liu et al.~\cite{liu2017} performed a shape optimisation of a $\Gamma$-shaped sound barrier, where the direct differentiation method (DDM) was employed to compute the sensitivity of the objective function with respect to CPs. Takahashi et al.~\cite{takahashi2019ewco}, which is a prior research of the current work, optimised periodic and layered structures in terms of the ultra-thin solar panels. They derived the shape derivatives with the adjoint variable method (AVM). Ummidivarapou et al.~\cite{ummidivarapu2020} introduced a teaching-learning-based optimisation algorithm, which is a gradient-free method, to design a acoustic horn. Similarly, Shaaban et al.~\cite{shaaban2020} performed a shape optimisation by exploiting the particle swarm optimisation (PSO) algorithm, which is gradient-free. This was extend to the axi-symmetric problem by the same authors~\cite{shaaban2020b}. The shape optimisation by Wang et al.~\cite{wang2020} is similar to Liu et al.~\cite{liu2017} but used the AVM instead of the DDM. On the other hand, the 3D acoustics was considered only by Chen et al.~\cite{chen2019}. They conducted a shape optimisation based on the DDM. Thus, their study can be regarded as a 3D version of \cite{liu2017}. They maximised the sound pressure of the surface of submarine or vase successfully.

Similarly to Chen et al.~\cite{chen2019}, the purpose of this study is to establish a shape optimisation system for 3D acoustic problems. In this system, a nonlinear optimisation algorithm integrates the corresponding IGBEM and AVM. These two ingredients were developed in the authors' previous research~\cite{takahashi2018jascome}. They proposed an accurate method to evaluate the singular and nearly-singular integrals associated with the isogeometric discretisation and, additionally, performed a shape-sensitivity analysis as an application.
The present work makes steady progress toward the shape optimisation with considering some optimisation algorithms which are implemented in two software Ipopt~\cite{ipopt} and NLopt~\cite{NLopt}. Those algorithm are compared with respect to their performances in some numerical examples.

The rest of this paper is organised as follows: Section~\ref{s:igbem} overviews an IGBEM for the 3D Helmholtz equation in terms of exterior homogeneous Neumann problems, which was constructed in our previous work~\cite{takahashi2018jascome}. Section~\ref{s:sopt} formulates the shape optimisation on the basis of the IGBEM and the adjoint variable method and describes the reduction of the problem to a nonlinear optimisation problem. Section~\ref{s:num} validates the proposed shape optimisation system through a numerical example and then demonstrates the capability of the system for complicated problems. Finally, Section~\ref{s:concl} concludes the  present study.

\section{Isogeometric BEM}\label{s:igbem}

We will overview the formulation of the IGBEM for the 3D Helmholtz equation, referring to our previous work~\cite{takahashi2018jascome}.

\subsection{Problem statement}\label{s:problem}

Let us consider a scattering problem of the time-harmonic acoustic wave in 3D. Specifically, we will solve the following exterior Neumann boundary value problem (BVP) in the infinite domain $\bbbr^3\setminus\overline{V}$: 
\begin{subequations}
  \begin{align}
     &\text{Governing equation}: &&\triangle u + k^2 u = 0 && \text{in $\bbbr^3\setminus\overline{V}$}, \label{eq:helm3d}\\
    &\text{Boundary condition}: &&\frac{\partial u}{\partial n}=0 && \text{on $S$},\label{eq:bc}\\
    &\text{Radiation condition}: &&u(\bm{x})\rightarrow\uin(\bm{x}) && \text{as $\abs{\bm{x}}\rightarrow\infty$},
  \end{align}%
  \label{eq:primary}%
\end{subequations}
where $u:\bm{x}\in\bbbr^3\to\bbbc$ denotes the total field or sound pressure, $\uin$ denotes a given incident field, $V$ denotes one or more acoustically-hard scatterers in $\bbbr^3$, $S$ denotes the boundary $\partial V$, $\bm{n}$ denotes the unit outward normal to $S$ and $k$ denotes the prescribed wavenumber.

\subsection{Boundary integral equation}\label{s:bie}

We will solve the BVP in (\ref{eq:primary}) with the following standard boundary integral equation (BIE):
\begin{eqnarray}
  C(\bm{x})u(\bm{x})+\int_S\frac{\partial G(\bm{x}-\bm{y})}{\partial n_y} u(\bm{y})\diff S_y = \uin(\bm{x})\quad\text{for $\bm{x}\in S$},
  \label{eq:bie}
\end{eqnarray}
where $G$ denotes the fundamental solution of the 3D Helmholtz equation, that is,
\begin{eqnarray}
  G(\bm{x}):=\frac{e^{\imath k |\bm{x}|}}{4\pi|\bm{x}|}.
  \label{eq:G}
\end{eqnarray}
Also, $C$ denotes the free term and is equal to $1/2$ if $S$ is smooth at $\bm{x}$. In this study, we utilise the equi-potential condition to yield
\begin{eqnarray}
  C(\bm{x})=1-\int_S\frac{\partial\Gamma(\bm{x}-\bm{y})}{\partial n_y}\diff S_y,
  \label{eq:C}
\end{eqnarray}
where $\Gamma(\bm{x}):=\frac{1}{4\pi\abs{\bm{x}}}$ denotes the fundamental solution for the Laplace equation in 3D. 

\subsection{Isogeometric analysis}\label{s:iga}

We will discretise the BIE in (\ref{eq:bie}) as well as the RHS of (\ref{eq:C}) under the concept of the isogeometric analysis (IGA). The IGA is a kind of isoparametric formulation that exploits the NURBS basis as both interpolation and shape functions mainly in the field of both boundary and finite element methods.

First, we express a given boundary $S$, which is supposed to consist of one or more closed surfaces, by using multiple NURBS surfaces. Each NURBS surface, say $\Pi$, is parameterised with two curve parameters $s$ and $t$, where the domain of $s$ and $t$ can be $[0,1]$ without the loss of generality. Then, we can express any point $\bm{y}$ on $\Pi$ as the tensor product of NURBS basis as follows:
\begin{eqnarray}
  \bm{y}(s,t)
  =\frac{\displaystyle\sum_{k=0}^{\ns-1}\sum_{l=0}^{\nt-1}w_{kl}N_k^\ps(s)N_l^\pt(t)\bm{C}_{kl}}{\displaystyle\sum_{k'=0}^{\ns-1}\sum_{l'=0}^{\nt-1}w_{k'l'}N^\ps_{k'}(s)N^\pt_{l'}(t)}
  =\sum_{k,l}\frac{w_{kl}\NN_{kl}(s,t)}{W(s,t)}\bm{C}_{kl},
  \label{eq:y}
\end{eqnarray}
where $N_k^p$ denotes the $k$-th B-spline function of degree $p$ and $w_{kl}$ and $\bm{C}_{kl}$ denote the $(k,l)$-th weight and control points, respectively, which should be determined according to the shape of $S$. Also, for the sake of simplicity, we denote the product $N_k^\ps(s)N_l^\pt(t)$ by $\NN_{kl}(s,t)$ and the summation in the denominator by $W(s,t)$.

The two series of knots, which are denoted by $\{s_i\}_{i=0}^{\ns+\ps}$ and $\{t_j\}_{j=0}^{\nt+\pt}$, are non-decreasing in general. To guarantee that the outer control points, i.e. the control points $\bm{C}_{kl}$ whose index $k$ or $l$ is either $0$ or the largest one (i.e. $\ns-1$ or $\nt-1$), locate on the perimeter $\partial\Pi$ of the NURBS surface, we use the clamped knots, i.e.
\begin{eqnarray*}
  s_i=
  \begin{cases}
    0 & i=0,\ldots,\ps\\
    \frac{i-\ps}{\ns-\ps} & i=\ps+1,\ldots,\ns-1\\
    1 & i=\ns,\ldots,\ns+\ps
  \end{cases}
\end{eqnarray*}
for $s$ and the same for $t$.

Similarly to the boundary point $\bm{y}$ in (\ref{eq:y}), we interpolate the boundary density $u$ on a surface $\Pi$ with the tensor product of the NURBS basis as follows:
\begin{eqnarray}
  u(s,t)=\sum_{k,l}\frac{w_{kl}\NN_{kl}(s,t)}{W(s,t)}u_{kl},
  \label{eq:u}
\end{eqnarray}
where coefficients $u_{kl}$ are the unknown variables to be determined from the BIE in (\ref{eq:bie}).

It should be noted that, since the knots are clamped, $u$ at a control point $\bm{C}_{kl}$ on the perimeter $\partial\Pi$ corresponds to $u_{kl}$ exactly; meanwhile, the other coefficients $u_{kl}$ do not generally correspond to $u$ at $\bm{C}_{kl}$.

The solution, that is, Dirichlet data $u$ on $S$ must be continuous across the intersecting line between two adjacent NURBS surfaces. This continuity-requirement can be satisfied by giving a unique unknown index, say $\nu$, to all the unknown coefficients associated with the underlying intersection. For example, let us consider the case that an outer control points $\bm{C}_{kl}$ on a NURBS surface $\Pi$ has the same position as an outer point $\bm{C}'_{k'l'}$ on another surface $\Pi'$, where we measure the geometrical distance of the two points $\bm{C}_{kl}$ and $\bm{C}'_{k'l'}$ to judge if they share the same position or not. Then, we give a global unknown index $\nu$ to the two points $\bm{C}_{kl}$ and $\bm{C}'_{k'l'}$ as well as the corresponding unknown coefficients $u_{kl}$ and $u_{k'l'}'$. As a result, we can obtain a certain number $N$ that represents the number of (global) unknowns over $S$. By using $N$ global unknowns and control points denoted by $u_\nu$ and $\bm{C}_\nu$, respectively, we no longer use the local indices (i.e. $kl$ and $k'l'$) and can express any point $\bm{y}$ and the boundary value $u$ as follows:
\begin{eqnarray}
  \bm{y}(s,t)=\sum_{\nu=1}^N R_\nu(s,t)\bm{C}_\nu,\quad u(s,t)=\sum_{\nu=1}^N R_\nu(s,t)u_\nu,
  \label{eq:y_u}
\end{eqnarray}
where $R_\nu$ corresponds to the basis $\frac{w_{kl} N_{kl}}{W}$ for a certain NURBS surface.

\subsection{Discretisation of the BIE}

By plugging (\ref{eq:y_u}) into the BIE in (\ref{eq:bie}), we can yield the following discretised BIE:
\begin{eqnarray}
  C(\bm{x}(\hat{s},\hat{t}))\sum_{\nu=1}^N R_\nu(\hat{s},\hat{t})u_\nu
  +\int_S \frac{\partial G}{\partial n_y}(\bm{x}(\hat{s},\hat{t}),\bm{y}(s,t))\sum_{\nu=1}^N R_\nu(s,t)\diff S_y u_\nu
  =\uin(\bm{x}(\hat{s},\hat{t})).
    \label{eq:bie2}
\end{eqnarray}
Here, a pair of parameters $(\hat{s},\hat{t})$ corresponds to a collocation point $\bm{x}$ on $S$ and each parameter is determined as the Greville abscissa~\cite{liu2017}. Similarly to the determination of the global $N$ unknowns ($u_\nu$), we regard the repeated collocation point on an intersection as a unique collocation point. As a result, we can determine $N$ distinct collocation points on $S$, which are enough to solve (\ref{eq:bie2}). In this study, we use the LU decomposition to solve $N$ unknowns ($u_\nu$) from a set of $N$ discretised BIEs of (\ref{eq:bie2}).

Once the unknowns are obtained, we can compute $u$ at any point: we may use (\ref{eq:u}) for any point on $S$, while we may exploit the integral representation for any point in $V$. In addition, we can compute the derivatives of $u$ on $S$ by differentiating the NURBS functions in (\ref{eq:u}) with respect to $s$ and/or $t$. This is useful to compute the shape derivative (sensitivity) because it usually consists of the derivative(s) of $u$ on a surface, as seen in (\ref{eq:sd}).

Regarding the boundary integrals in (\ref{eq:bie2}), we apply the Lachat's method to the singular integrals and the hierarchical subdivision technique to the singular- and nearly-singular-integrals. The details are described in our previous paper~\cite{takahashi2018jascome}.

\subsection{Knot insertion}\label{s:ki}

As we will mention in Section~\ref{s:sopt}, we will optimise the shape of $S$ via the control points $\bm{C}_\nu$. If the number of control points involved in a target $S$ is large, the convergence of the optimisation would be slow. So, one may construct a surface with a small number $N$ of control points. However, this can lead to a low accurate solution in the IGA because the number of unknowns (degrees of freedom) is also $N$; recall (\ref{eq:y_u}). To resolve this issue, which is common in the IGA~\cite{hughes2005}, we may resort to the knot insertion, by which control points can be added to anyplace without changing the shape of $S$. This technique is used when we analyse the BVP in (\ref{eq:primary}) as well as the adjoint one in (\ref{eq:adjoint}), which will be mentioned in Section~\ref{s:statement}.
\def\mathcalS{\mathcal{S}}
\def\mathcalJ{\mathcal{J}}
\section{Gradient-based shape optimisation}\label{s:sopt}
We will construct our shape optimisation method based on the IGBEM, adjoint variable method and nonlinear optimisation method. The present framework is a direct extension of the 2D case investigated in our previous paper~\cite{takahashi2019ewco}\footnote{The corresponding software is available at \texttt{https://sourceforge.net/projects/igbemsopt/}.}.

\subsection{Problem statement and shape derivative}\label{s:statement}

The present shape optimisation problem is to maximise or minimise a prescribed objective function $\mathcalJ$ by changing the surface of scatterer(s) $V$, i.e. the boundary $S$. Specifically, we define $\mathcalJ$ as the summation of the sound pressure $u$ at $M$ observation points $\{\bm{z}_i\}_{i=1}^M$, that is,
\begin{eqnarray}
  \mathcalJ(u;S):=\sum_{m=1}^M\frac{\abs{u(\bm{z}_m)}^2}{2},
  \label{eq:J}
\end{eqnarray}
where $u$ is supposed to be a solution of the BVP in (\ref{eq:primary}) or the \textit{primary} problem in the context of the adjoint variable method.

To define the shape derivative (sensitivity), denoted by $\mathcalS$, of $\mathcalJ$ in (\ref{eq:J}), we slightly move every point $\bm{y}$ on $S$ by $\epsilon\bm{V}(\bm{y})$, where $\epsilon$ is an infinitesimally small number and $\bm{V}$ denotes the direction to move. Correspondingly, the boundary $S$ and the field $u$ are perturbated to $\tilde{S}$ and $\tilde{u}$, respectively. Then, $\mathcalS$ is defined as the coefficient of the term $\epsilon$ which is obtained by expanding the perturbated objective function $\mathcalJ(\tilde{u}; \tilde{S})$ with respect to $\epsilon$. Therefore, we have
\begin{eqnarray}
  \mathcalJ(\tilde{u};\tilde{S})=\mathcalJ(u;S)+\epsilon\mathcalS(u;S)+O(\epsilon^2).
  \label{s:J_expand}
\end{eqnarray}
Here, as well-known (see \cite{feijoo2003} for example), $\mathcalS$ can be derived as follows:
\begin{eqnarray}
  \mathcalS(u;S)=\Re\int_S\left(k^2\lambda^* u-\nabla \lambda^* \cdot \nabla u\right)\bm{V}\cdot\bm{n}\ \diff S,
  \label{eq:sd}
\end{eqnarray}
where $()^*$ denotes the complex conjugate and the adjoint field $\lambda$ is the solution of the following adjoint problem:
\begin{subequations}
  \begin{align}
    &\text{Governing equation}: &&\triangle \lambda(\bm{x}) + k^2 \lambda(\bm{x}) = -\sum_m u(\bm{x})\delta(\bm{x}-\bm{z}_m) && \text{in $V$},\label{eq:helm_adj}\\
    &\text{Boundary condition}: &&\frac{\partial\lambda}{\partial n}=0 && \text{on $S$},\label{eq:bc_adj}\\
    &\text{Radiation condition}: &&\lambda(\bm{x})\rightarrow 0 && \text{as $\abs{\bm{x}}\rightarrow\infty$}.\label{eq:radiation_adj}%
  \end{align}%
  \label{eq:adjoint}%
\end{subequations}
We can also solve the adjoint problem with the IGBEM mentioned in Section~\ref{s:igbem}. To this end, we may replace the incident field $\uin(\bm{x})$ in (\ref{eq:bie}) with the following term:
\begin{eqnarray*}
  -\int_V G(\bm{x}-\bm{y})\left(-\sum_m u(\bm{y})\delta(\bm{y}-\bm{z}_m)\right)\diff V_y
  = \sum_m G(\bm{x}-\bm{z}_m)u(\bm{z}_m),
\end{eqnarray*}
where $G$ is the fundamental solution given in (\ref{eq:G}).

\subsection{Discretisation of the shape derivative}

In the numerical analysis, the infinitesimal deformation (perturbation) $\epsilon\bm{V}$ must be finite. When we denote the point of an arbitrary point $\bm{y}$, which is expressed as (\ref{eq:y_u}) in the IGBEM, on the surface $\tilde{S}$ by $\tilde{\bm{y}}$, we can approximate $\epsilon\bm{V}(\bm{y})$ as follows:
\begin{eqnarray*}
  \epsilon\bm{V}(\bm{y})\approx\tilde{\bm{y}}(s,t)-\bm{y}(s,t)
  =\sum_{\nu=0}^{N-1} R_\nu(s,t)\delta\bm{C}_\nu,
\end{eqnarray*}
where $\delta\bm{C}_{\nu}$ denotes the variation of the control point $\bm{C}_\nu$, i.e.
\begin{eqnarray*}
  \delta\bm{C}_{\nu}:=\tilde{\bm{C}_{\nu}}-\bm{C}_{\nu}.
\end{eqnarray*}
Then, $\mathcalJ$ in (\ref{s:J_expand}) can be discretised as follows:
\begin{eqnarray}
  \mathcalJ(\tilde{u};\tilde{S})\approx \mathcalJ(u;S)+\sum_{\nu=0}^{N-1}\bm{s}_\nu(u;S)\cdot\delta\bm{C}_\nu+O(\epsilon^2),
  \label{eq:tildeJ}
\end{eqnarray}
where
\begin{eqnarray*}
  \bm{s}_\nu(u;S):=\Re\int_S\left( k^2\lambda^* u-\nabla \lambda^* \cdot \nabla u
\right) R_\nu \bm{n}\diff S.
  \label{eq:sd_disc}
\end{eqnarray*}
The vector $\bm{s}_\nu$ stands for the sensitivity of $\mathcalJ$ with respect to the control points $\bm{C}_\nu$, which are considered as the design variables in this study.

Because of $\partial u/\partial n=0$ and $\partial\lambda/\partial n=0$ due to the boundary conditions in (\ref{eq:bc}) and (\ref{eq:bc_adj}), respectively, the gradients of $u$ and $\lambda^*$ in (\ref{eq:sd_disc}) can be expressed as follows:
\begin{eqnarray*}
  \nabla u = \frac{1}{J}\left(\pp{u}{s}\bm{t}\times\bm{n}+\pp{u}{t}\bm{n}\times\bm{s}\right),\quad
  \nabla \lambda^* = \frac{1}{J}\left(\pp{\lambda^*}{s}\bm{t}\times\bm{n}+\pp{\lambda^*}{t}\bm{n}\times\bm{s}\right),
\end{eqnarray*}
where $\bm{s}:=\pp{\bm{y}}{s}$ and $\bm{t}:=\pp{\bm{y}}{t}$ denote the tangential vectors along $s$ and $t$ coordinates, respectively, and $J:=\norm{\bm{s}\times\bm{t}}$ is the Jacobian. It should be emphasised that the tangential derivatives $\pp{u}{s}$ and $\pp{u}{t}$ can be computed readily by differentiating the NURBS basis $R_\nu$. In addition, the gradients are continuous over the surface $S$ (except for the intersections among NURBS surfaces in general) if the degrees $\ps$ and $\pt$ of the NURBS basis are two or more.

Since there is no singularity in the integral in (\ref{eq:sd_disc}), we may evaluate the integral with the Gauss-Legendre quadrature formula.

\subsection{Reduction to nonlinear optimisation problem}\label{s:nonlinear}

The optimisation problem stated in Section~\ref{s:statement} forms a nonlinear optimisation problem. In general, the problem is to minimise the prescribed objective function $f:\bbbr^n\rightarrow\bbbr$ with respect to $n$ design variables $x\in\bbbr^n$ under $m$ inequality-constraints $g:\bbbr^n\rightarrow\bbbr^m$, i.e. 
\begin{eqnarray}
  g_{\rm L} \leq g(x)\leq g_{\rm U},
\end{eqnarray}
where $g_{\rm L,U}\in\bbbr^m$ denote the bounds of $g$. The design variables $x$ are usually bounded as
\begin{eqnarray}
  x_{\rm L} \leq x\leq x_{\rm U},
\end{eqnarray}
where $x_{\rm L,U}\in\bbbr^m$ denote the bounds of $x$. Optionally, the gradient $\nabla f$ and Hessian of $f$ is considered if they are readily computed.

In the present shape optimisation problem, we choose the $N$ control points $\bm{C}_\nu$ (where $\nu=0,\ldots,N-1$) as the design variables. Then, we may regard our objective function $\mathcalJ$, design variables $\bm{C}_\nu$ and their gradients $\bm{s}_\nu$ in (\ref{eq:sd_disc}) as $f$,  $(\bm{C}_0,\ldots,\bm{C}_{N-1})^{\rm T}\in\bbbr^{3N}$ and $(\bm{s}_0,\ldots,\bm{s}_{N-1})^{\rm T}\in\bbbr^{3N}$, respectively, where $3N$ corresponds to the number $n$ of design variables. In this study, we do not consider the Hessian of $\mathcalJ$.

We utilise a primal-dual interior-point method with line searches based on Filter methods, which is implemented in Ipopt~\cite{ipopt,ipopt_wiki} and will be called IP hereafter. In the previous research~\cite{takahashi2019ewco}, the IP sometimes required many number of backtracking line-search steps and, thus, many number of evaluating $f$ as well as $\nabla f$. As a result, the computational time was sometimes enormous.

Hence, we consider different gradient-based optimisation methods. To this end, we exploit the software NLopt~\cite{NLopt}, which contains many optimisation methods. We use the MMA (method of moving asymptotes~\cite{svanberg2002class}) and SLSQP (sequential least-squares quadratic programming~\cite{kraft1994}) because they are general-purpose in the sense that they can handle nonlinear inequality-constraints.

\def\mathcalS{\mathcal{S}}
\def\mathcalJ{\mathcal{J}}

\section{Numerical examples}\label{s:num}

This section will provide some numerical examples with our shape optimisation software. In Section~\ref{s:validate}, we will validate the software through an optimisation problem, which can be analysed exactly. The problem is actually a parametric optimisation problem, but we can test the entire of our software for the shape optimisation problem. After the validation, we will show the capability of the software through some more complicated examples in Sections~\ref{s:reflector}--\ref{s:bending}.

\subsection{Verification}\label{s:validate}

\nprounddigits{5}

\subsubsection{Problem configuration}

\begin{figure}[H]
  \centering
  \includegraphics[width=.35\textwidth]{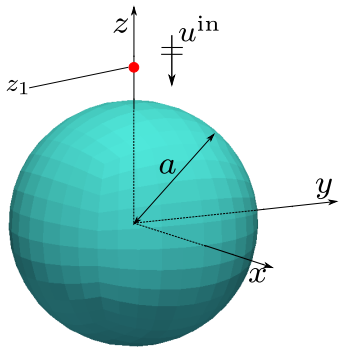}
  \caption{Problem setting.}
  \label{fig:test}
\end{figure}

Let us consider a parametric optimisation problem. Specifically, we will find the radius, denoted by $a$, of a spherical scatterer (centred at the origin) so that the radius can maximises an objective function $\mathcalJ$ in (\ref{eq:J}), where $M=1$ and the corresponding observation point $\bm{z}_1$ is chosen as $(0, 0, 8.5)^{\rm T}$ (Figure~\ref{fig:test}): namely,
\begin{eqnarray}
  \mathcalJ(u;S):=\frac{\abs{u(\bm{z}_1)}^2}{2}\quad\text{where $\bm{z}_1=(0, 0, 8.5)^{\rm T}$.}
  \label{eq:J_verification}
\end{eqnarray}
We give a planewave incident field $\uin(\bm{x})=\mathrm{e}^{-\imath kz}$, which propagates in the $-z$ direction, where the wavenumber $k$ is given as one.

Following the reference~\cite{cobb1988}, we create the surface $S$ of the spherical scatter with six NURBS surfaces. Each surface is constructed with $5\times 5$ control points and the tensor product of the B-spline functions of degree 4, i.e. $\ns=\nt=5$ and $\ps=\pt=4$ . Correspondingly, the number $N$ of the (unique) control points is 98. As noted in Section~\ref{s:ki}, we can increase $N$ to improve the resolution of the boundary element solution by the knot insertion. We consider three cases of $N$, i.e. $N=866$, $2402$ and $3458$.

Under the present configuration, the sound pressure $u$ at the observation point $\bm{z}_1$ can be written as a function of the radius $a$~\cite{bowman1987}, that is,
\begin{eqnarray}
  u(\bm{z}_1) = \sum_{n=0}^{\infty} i^n (2n + 1)(j_n(kr) - A'_n {h^{(1)}_n}'(kr))P_n (\cos(\theta-\pi)),
  \label{eq:uanal}
\end{eqnarray}
where the spherical coordinates $\theta$ and $r$ corresponding to $\bm{z}_1$ are $0$ and $8.5$, respectively.
 Also, $j_n$, $h^{(1)}_n$ and $P_n$ denote the spherical Bessel function of degree $n$,  the spherical Hankel function of the first kind and degree $n$ and the Legendre polynomial of degree $n$, respectively. In addition, the coefficient $A'_n$ is defined as
\begin{eqnarray*}
  A'_n := \frac{j'_n (ka)}{{h^{(1)}_n}'(ka)},
\end{eqnarray*}
where the prime represents the differentiation with respect to $a$. It should be noted that (\ref{eq:uanal}) is valid when the observation points $\bm{z}_1$ is in the outside of the sphere, that is, $a\in(0,8.5)$.

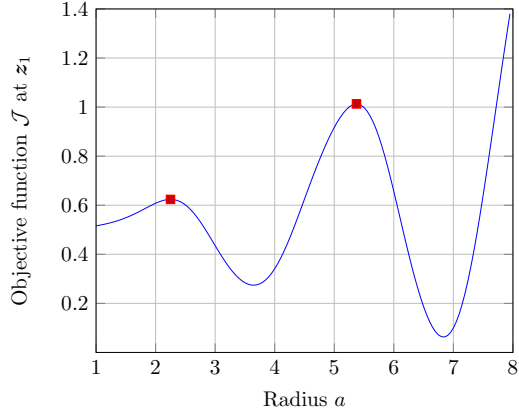
\begin{figure}[H]
  \centering
  \begin{tikzpicture}[scale=0.8]
  \begin{axis}[
      xlabel={Radius $a$},
      ylabel style={align=center},
      ylabel={Objective function $\mathcalJ$ at $\bm{z}_1$},
      xmin=1.0, xmax=8.0,
      ymin=0.0, ymax=1.4, ytick={0.2,0.4,...,2.0},
      xmajorgrids=true,
      ymajorgrids=true,
    ]
    \addplot table [sharp plot, mark=none, x=a, y=ana] {data/J.txt};
    \addplot table [only marks] {
      2.24991750058038420e+00 6.23794337070834093e-01
      5.37318093130602481e+00 1.01260319466883364e+00
    };
  \end{axis}
\end{tikzpicture}
  \caption{The value of the objective function $\mathcalJ$ in (\ref{eq:J_verification}) as a function of the radius $a$ of the spherical scatterer at the origin. The two points represent the local maxima in (\ref{eq:maxima}) when $a$ is restricted to $[1,7]$.}
  \label{fig:J_a}
\end{figure}

From (\ref{eq:J_verification}) and (\ref{eq:uanal}), we can plot $\mathcalJ$ against $a\in[1,8]$ as in Figure~\ref{fig:J_a}. When we restrict the lower and upper bounds of the design variable $a$ to $1$ and $7$, respectively, we can have two local maxima, i.e.
\begin{eqnarray}
  (a,J)=(\numprint{2.24991750058038420e+00},\ \numprint{0.623794337070834093}),\quad
  (\numprint{5.37318093130602481e+00},\ \numprint{1.01260319466883364e+00}),
  \label{eq:maxima}
\end{eqnarray}
which are computed by applying the Brent minimisation algorithm, which is implemented in GNU Scientific Library~\cite{GSL}, to (\ref{eq:uanal}). If we can arrive at one of these local maxima from a certain initial radius, denoted by $a_0$, we can validate our optimisation software. In what follows, we consider two initial values, i.e. $a_0=3$ and $4$.

The design variable of this problem is the radius $a$ only, while our shape optimisation method treats all the control points $\bm{C}_\nu$ as the design variables (recall Section~\ref{s:nonlinear}). To fill this gap, we modify a fraction of the computer program by considering the relationship between the radius $a$ and control points $\bm{C}_\nu$. 
Specifically, when the radius is updated from $a$ to $\tilde{a}$, all the control points must be scaled by $\tilde{a}/a$ so that the surface $S$ can preserve its spherical shape. Therefore, the control points $\bm{C}_\nu$ must be updated to $\tilde{\bm{C}}_\nu$ so that
\begin{eqnarray}
  \tilde{\bm{C}}_\nu = \frac{\tilde{a}}{a} \bm{C}_\nu
  \label{eq:tildeC}
\end{eqnarray}
holds. Plugging this and the variation of the radius, i.e. $\delta a:=\tilde{a}-a$, into (\ref{eq:tildeJ}), we have
\begin{eqnarray}
  \mathcalJ(\tilde{u};\tilde{S}) = \mathcalJ(u;S)+\sum_{\nu=1}^N \frac{\bm{s}_\nu(u;S)\cdot\bm{C}_\nu}{a}\delta a+O(\epsilon^2).
\end{eqnarray}
Clearly, the shape derivative of $\mathcalJ$ with respect to the radius $a$, i.e. $\pp{\mathcalJ}{a}$, is obtained as
\begin{eqnarray}
  \pp{\mathcalJ}{a}=\sum_{\nu=1}^N \frac{\bm{s}_\nu(u;S)\cdot\bm{C}_\nu}{a}.
  \label{eq:djda}
\end{eqnarray}

Hence, when the new (perturbated) radius $\tilde{a}$ is determined by an optimisation algorithm, we first update the control points $\bm{C}_\nu$ to $\tilde{\bm{C}}_\nu$ according to (\ref{eq:tildeC}) and, then, $\pp{\mathcalJ}{a}$ in (\ref{eq:djda}). This procedure is added to the user-defined routine to compute $\mathcalJ$ and its gradient.

In this analysis, we let the convergence tolerance of IP, i.e. the parameter \texttt{tol} of Ipopt~\cite[Page 69]{ipopt} be $10^{-3}$.
In regard to both MMA and SLSQP, we let the relative tolerance, which corresponds to the parameter \texttt{ftol\_rel} of NLopt, be $10^{-3}$.

\subsubsection{Results and discussions}

Table~\ref{tab:sphere} shows the computed optimal radius $a$ and the corresponding value of $\mathcalJ$ for the two initial radii $a_0$ and three numbers $N$ of control points.
We can observe that every solver could achieve one of the local maxima.
There is no clear difference due to $N$, which means that the discretisation error of the IGBEM is almost negligible with the smallest $N$.

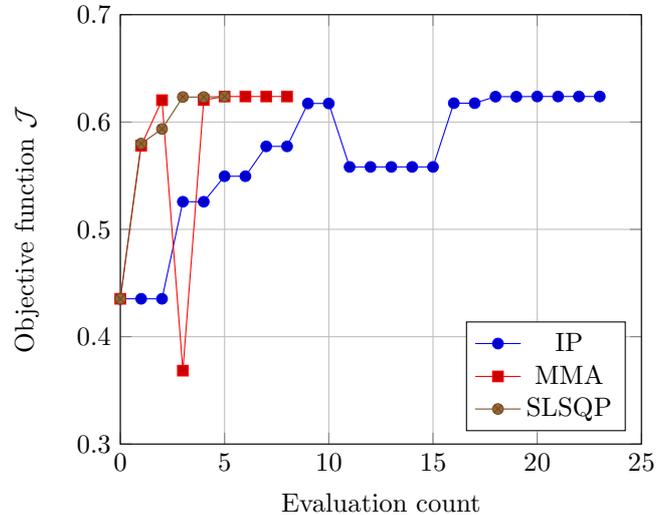
\begin{figure}[H]
  \centering
  \begin{tikzpicture}[scale=1.0]
  \begin{axis}[
    xlabel={Evaluation count},
    ylabel style={align=center},
    ylabel={Objective function $\mathcalJ$},
    xmin=0, xmax=25, xtick={0,5,...,25},
    ymin=0.3, ymax=0.7, ytick={0.3,0.4,...,0.7},
    legend pos=south east,
    xmajorgrids=true,
    ymajorgrids=true,
    ]
    \addplot table [sharp plot, x=count, y=f] {data/nki8_a3.0_ipopt_CHECK_COUNT.txt};
    \addlegendentry{IP}
    \addplot table [sharp plot, x=count, y=f] {data/nki8_a3.0_nlopt24_CHECK_COUNT.txt};
    \addlegendentry{MMA}
    \addplot table [sharp plot, x=count, y=f] {data/nki8_a3.0_nlopt40_CHECK_COUNT.txt};
    \addlegendentry{SLSQP}
  \end{axis}
\end{tikzpicture}
  \caption{History of the value of the objective function $\mathcalJ$ in (\ref{eq:J_verification}) against the evaluation count in the case of the initial radius $a_0=3$ and $N=866$.}
  \label{fig:test_history}
\end{figure}

The columns of ``Eva.'' in Table~\ref{tab:sphere} shows the number of evaluating $\mathcalJ$ (and most likely its gradient at the same time) until convergence. In every combination of $a_0$ and $N$, the SLSQP required less evaluation counts than the others. Figure~\ref{fig:test_history} plots the value of $\mathcalJ$ against the evaluation count in the case of $a_0=3$ and $N=866$.

The present results indicate that our formulation and its numerical implementation are valid.

\begin{table}[H]
  \centering
  \caption{Result of the optimised radius $a$ and the corresponding value of the objective function $\mathcalJ$ in Section~\ref{s:validate}. Here, the heading ``Eva.'' stands for the number of evaluating $\mathcalJ$.}
  \label{tab:sphere}
  \scriptsize
  $a_0=3$\\
  \begin{tabular}{|c|c|c|c|c|c|c|c|c|c|}
    \hline
    Algo. & \multicolumn{3}{c|}{IP} &  \multicolumn{3}{c|}{MMA} &  \multicolumn{3}{c|}{SLSQP}\\
    \hline
    $N$ & $a$ & $\mathcalJ$ & Eva.  & $a$ & $\mathcalJ$ & Eva.  & $a$ & $\mathcalJ$ & Eva. \\
    \hline
    866 & \numprint{2.24998615900516929e+00} & \numprint{0.623794283354436985} & 23 & \numprint{2.24978909861472687e+00} & \numprint{0.623794279670119023} & 8 & \numprint{2.24274042762932035e+00} & \numprint{0.623778143104502414} & 5\\
    2402 & \numprint{2.24998081149424500e+00} & \numprint{0.623794335188337157} & 23 & \numprint{2.24978606250787783e+00} & \numprint{0.623794331009669412} & 8 & \numprint{2.24273797877532655e+00} & \numprint{0.623778182868081066} & 5\\
    3458 & \numprint{2.24998289952717334e+00} & \numprint{0.623794335575214687} & 23 & \numprint{2.24978730438999408e+00} & \numprint{0.623794331583398809} & 8 & \numprint{2.24273916078656965e+00} & \numprint{0.623778188638723918} & 5\\
    \hline
  \end{tabular}\\[\baselineskip]
  $a_0=4$\\
  \begin{tabular}{|c|c|c|c|c|c|c|c|c|c|}
    \hline
    Algo. & \multicolumn{3}{c|}{IP} &  \multicolumn{3}{c|}{MMA} &  \multicolumn{3}{c|}{SLSQP}\\
    \hline
    $N$ & $a$ & $\mathcalJ$ & Eva.  & $a$ & $\mathcalJ$ & Eva.  & $a$ & $\mathcalJ$ & Eva. \\
    \hline
    866 & \numprint{5.37302035529867261e+00} & \numprint{1.01254399549004437e+00} & 13 & \numprint{5.37327557417074431e+00} & \numprint{1.01254397411240737e+00} & 10 & \numprint{5.37319113655478908e+00} & \numprint{1.01254399352844837e+00} & 9\\
    2402 & \numprint{5.37295246736562149e+00} & \numprint{1.01260217818399623e+00} & 13 & \numprint{5.37315952650911122e+00} & \numprint{1.01260222197983096e+00} & 16 & \numprint{5.37312574673180787e+00} & \numprint{1.01260221984283727e+00} & 9\\
    3458 & \numprint{5.37298806847122368e+00} & \numprint{1.01260290777094331e+00} & 13 & \numprint{5.37320910135166496e+00} & \numprint{1.01260293879015717e+00} & 10 & \numprint{5.37315980955157002e+00} & \numprint{1.01260293911576449e+00} & 9\\
    \hline
  \end{tabular}
\end{table}

\subsection{Example 1: Reflector}\label{s:reflector}

\begin{figure}[H]
  \centering
  \includegraphics[width=.7\textwidth]{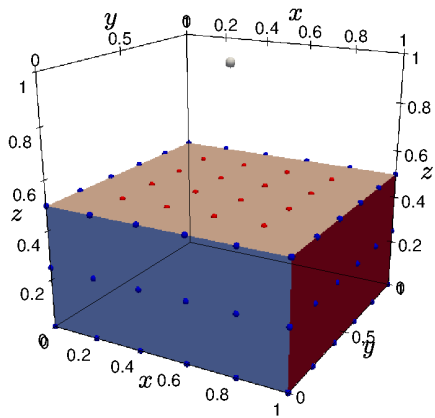}
  \caption{The initial shape of the reflector model in Section~\ref{s:reflector}. The point in grey represents the observation point $\bm{z}_1=(0.5,0.5,1.0)^{\rm T}$. The red and points represent the CPs that are designed, while the blue ones represent the fixed CPs.}
  \label{fig:app12a_config}
\end{figure}

To demonstrate the capability of our shape optimisation framework, we will begin with a simple model of a cuboid ``reflector'', whose dimensions are $1\times 1\times 0.5$, as shown in Figure~\ref{fig:app12a_config}. 
Regarding a planewave incident field $\uin$ with the wavenumber of $k=3$ that propagates in the $-z$ direction, i.e. $\uin(\bm{x})=e^{-\imath k z}$, we try to maximise the objective function $\mathcalJ$ in (\ref{eq:J}), where a single observation point $\bm{z}_1=(0.5,0.5,1.0)^{\rm T}$ in the illuminated side is considered.

The surface of the reflector consists of six NURBS surfaces (viz., top, bottom, left, right, front and back) with the NURBS functions of degree 2, i.e. $\ps=\pt=2$. They are shown in different colours in Figure~\ref{fig:app12a_config}. The number $\ns$, which denotes the number of control points along the local coordinate $s$, is given as $6$, $6$ and $3$ if $s$ is parallel to $x$-, $y$- and $z$-axis, respectively, at the initial configuration. Similarly, we determine the value of $\nt$ every NURBS surface. Each NURBS surface is clamped on its perimeter, as mentioned in Section~\ref{s:iga}. In this case, the number $N$ of unique CPs is 92, which are shown as points (in red or blue) on the surfaces in Figure~\ref{fig:app12a_config}. The number of CPs is increased to 548 by the knot insertion when we perform the isogeometric boundary element analysis.

In this example, we design only $4\times 4$ CPs, which are coloured in red in Figure~\ref{fig:app12a_config}, on the top surface excluding the perimeter. In addition, we allow each target CP to move vertically at most 0.3, which guarantees that any target CP never touches with the others. Thus, the number of design variables is 16. To be specific, we first regard all the coordinates of the CPs as the design variables and, then, set the initial coordinate to both lower and upper bounds for each coordinate that is not optimised.

\begin{figure}[H]
  \centering
  \begin{tikzpicture}[scale=1.0]
  \begin{axis}[
    xlabel={Evaluation count},
    ylabel style={align=center},
    ylabel={Objective function $\mathcalJ$},
    xmin=0, xmax=25, xtick={0,5,...,25},
    ymin=0.6, ymax=1.4, ytick={0.6,0.7,...,1.4},
    xmajorgrids=true,
    ymajorgrids=true,
    legend pos=south east,
    ]
    \addplot table [sharp plot, x=count, y=f] {data/ipopt_tol1e-3_CHECK_COUNT.txt};
    \addlegendentry{IP}
    \addplot table [sharp plot, x=count, y=f] {data/nlopt24_CHECK_COUNT.txt};
    \addlegendentry{MMA}
    \addplot table [sharp plot, x=count, y=f] {data/nlopt40_CHECK_COUNT.txt};
    \addlegendentry{SLSQP}
  \end{axis}
\end{tikzpicture}
  \caption{History of the value of the objective function $\mathcalJ$ in the reflector model (Section~\ref{s:reflector}).}
  \label{fig:app12a_history}
\end{figure}
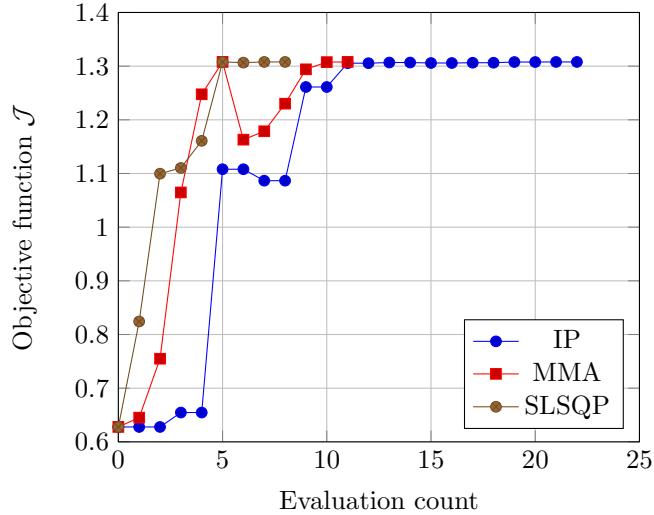

Figure~\ref{fig:app12a_history} compares the history of the value of $\mathcalJ$ for the three optimisation algorithms. All the algorithms converged to almost the same solution. Similarly to the previous example, the SLSQP required the least number of evaluations until convergence.

\begin{figure}[H]
  \centering
  \begin{tabular}{cc}
    \includegraphics[width=.4\textwidth]{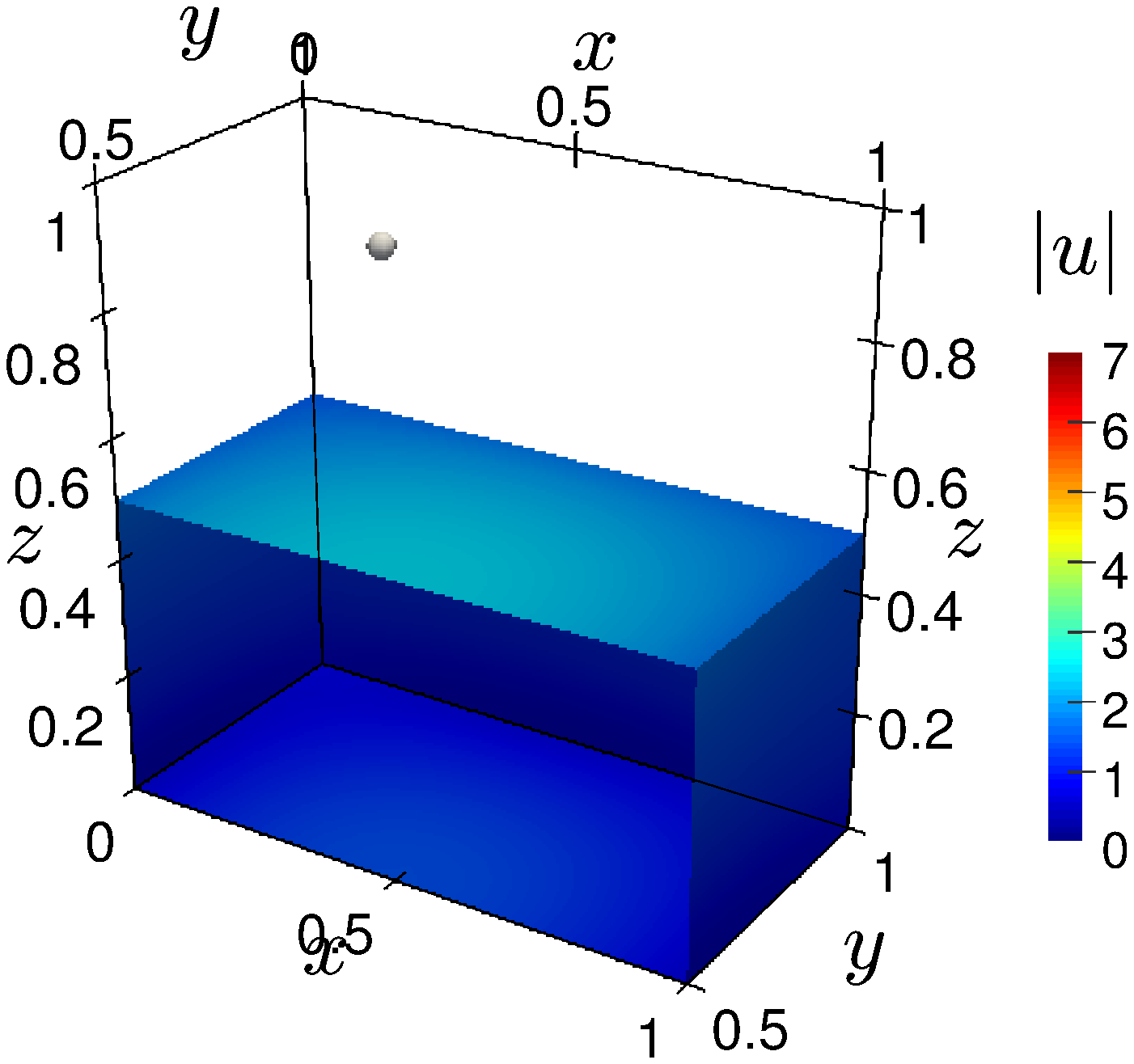}
    &\includegraphics[width=.4\textwidth]{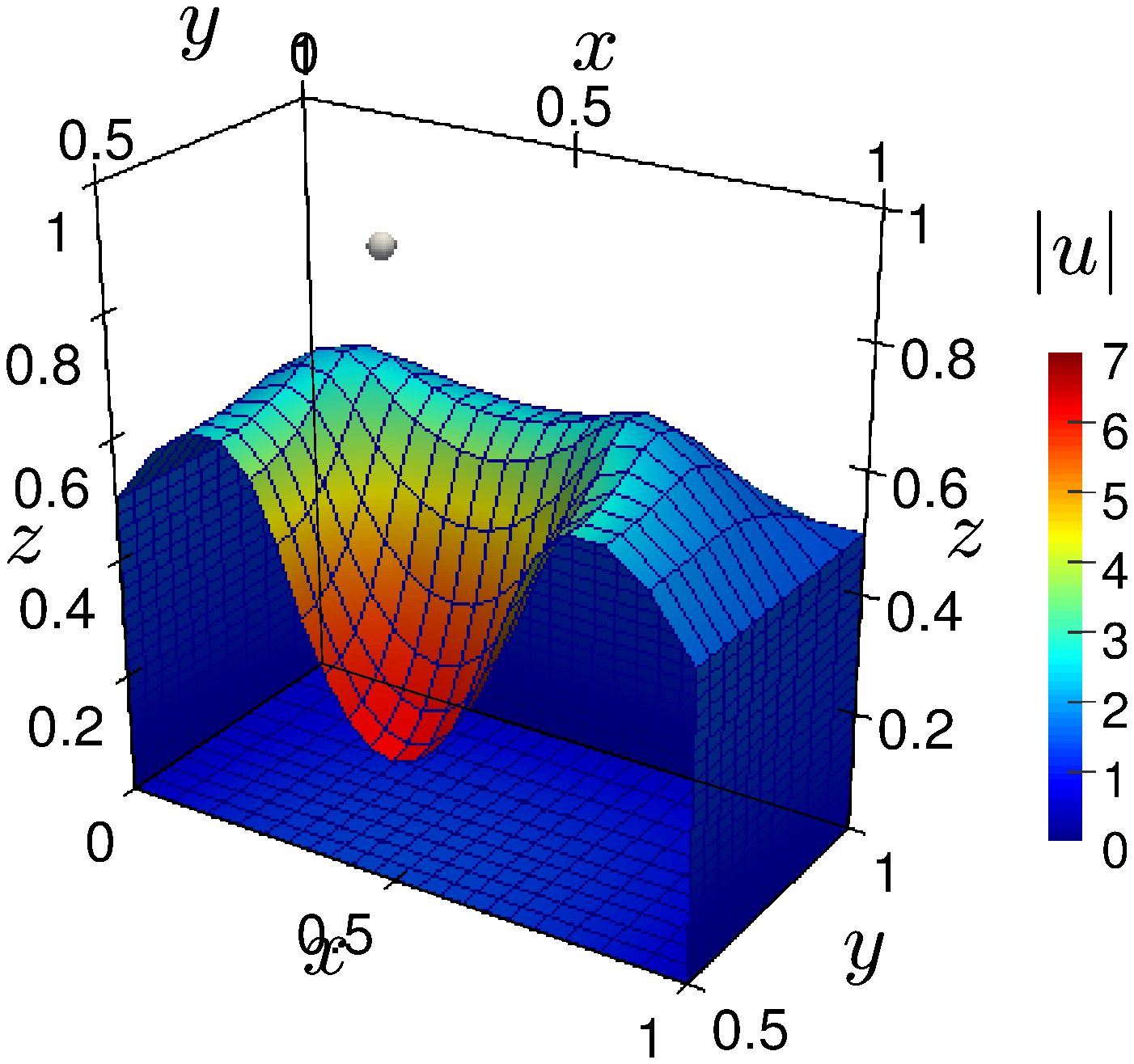}\\
    Initial shape & Optimised shape
  \end{tabular}
  \caption{Distribution of $|u|$ on the surface in the reflector model (Section~\ref{s:reflector}). The point (coloured in grey) is the observation point $\bm{z}_1=(0.5,0.5,1.0)^{\rm T}$.}
  \label{fig:app12a_u_abs_final}
\end{figure}

\begin{figure}[H]
  \centering
  \includegraphics[width=.8\textwidth]{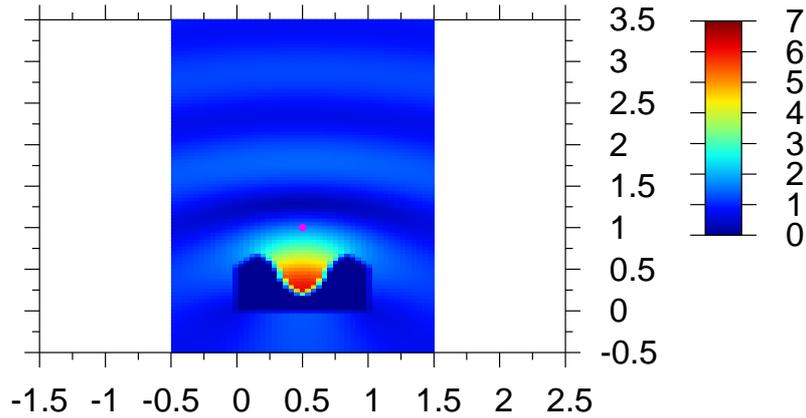}
  \caption{Distribution of $|u|$ on the plane of $y=0.5$ in the optimised shape of the reflector model (Section~\ref{s:reflector}).}
  \label{fig:app12a_uin_final}
\end{figure}

Figure~\ref{fig:app12a_u_abs_final} shows the distribution of the absolute value of sound pressure, i.e.  $|u|$, on the boundary at both initial and optimised shapes. In addition, Figure~\ref{fig:app12a_uin_final} shows the distribution of $|u|$ on the middle cross section, i.e. $y=0.5$; draw range was selected as $-0.5\le x\le 1.5$ and $-0.5\le z\le 3.5$. The peak of $|u|$ is not on the observation point $\bm{z}_1$, but the created shape is reasonable in the sense that it looks like a parabolic antenna. It should be noted that the results in Figures~\ref{fig:app12a_u_abs_final} and \ref{fig:app12a_uin_final} are of the SLSQP but almost the same results were obtained by both IP and MMA.

By considering the evaluation counts in Figure~\ref{fig:app12a_history} as well as Figure~\ref{fig:test_history}, we will use the SLSQP only in the following examples.

\subsection{Example 2: Resonator}\label{s:resonator}

\begin{figure}[H]
  \centering
  \begin{tabular}{cc}
    \includegraphics[width=.4\textwidth]{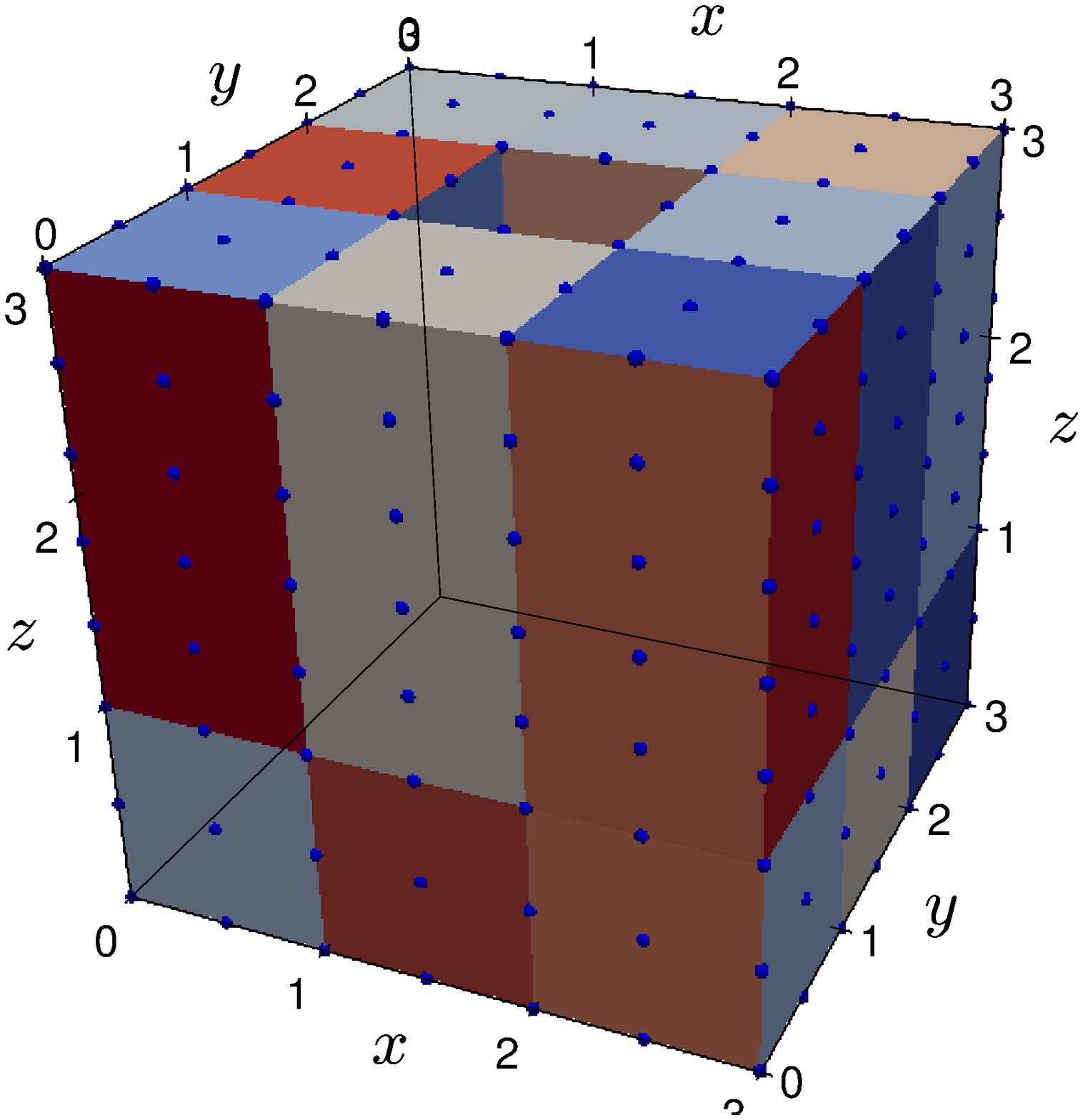}
    &\includegraphics[width=.4\textwidth]{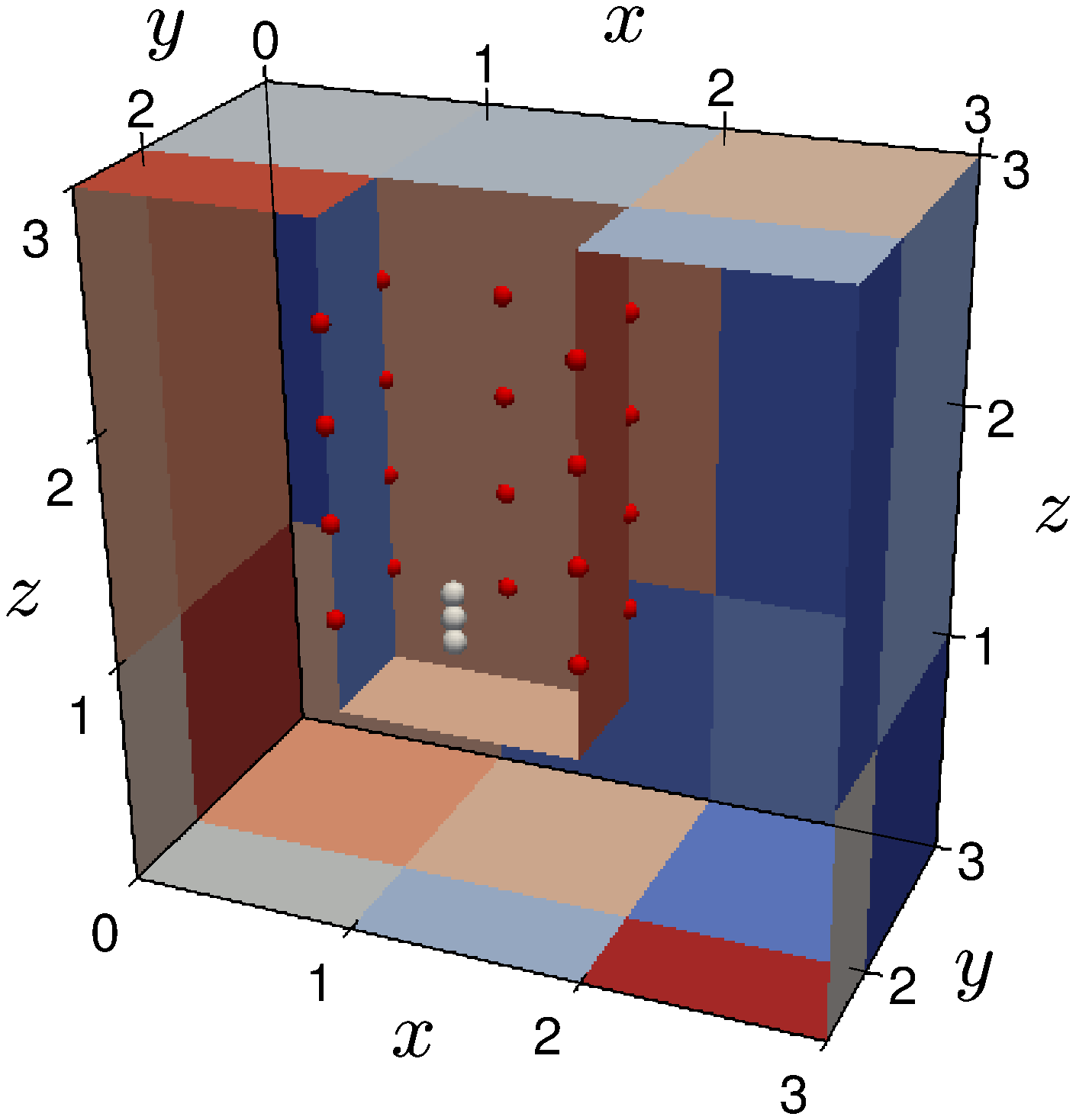}\\
    Entire & Posterior half
  \end{tabular}
  \caption{The initial shape of the resonator model in Section~\ref{s:resonator}. The points in grey represent the three observation point $\bm{z}_1$, $\bm{z}_2$ and $\bm{z}_3$ in the hollow. The red and points represent the CPs that are designed, while the blue ones represent the fixed CPs; the CPs on the left, back and bottom surfaces are also fixed.}
  \label{fig:app7f_config}
\end{figure}

As the second example, we attempt to catch a sound with a bowl. Specifically, as illustrated in Figure~\ref{fig:app7f_config} (left), we consider a cubic scatterer (whose dimensions are $3\times 3\times 3$) with a hollow ($1\times 1\times 2$). Then, we optimise the shape of the hollow so that we can increase the sound pressure inside it. We give two types of the incident fields, i.e. $\uin(\bm{x})=e^{-\imath k z}$ and $e^{-\imath k x}$, which represent the planewave propagating in the $-z$ and $-x$ direction, respectively. Here, $k=3$ is supposed.

Regarding the objective function $\mathcalJ$ in (\ref{eq:J}) to be maximised, we consider three observation points aligned vertically in the hollow, i.e. $\bm{z}_1=(1.5,1.5,1.4)$, $\bm{z}_2=(1.5,1.5,1.5)$ and $\bm{z}_3=(1.5,1.5,1.6)$, which are shown as three points (in grey) Figure~\ref{fig:app7f_config} (right).

The bowl is modelled with $46$ NURBS surfaces of degree 2, which are distinguished by different colours as in Figure~\ref{fig:app7f_config} (left). The whole surface of the bowl includes $282$ CPs, which are displayed as the points (in blue) in the same figure.  In addition, we increase the number $N$ of CPs from 282 to 1314 by the knot insertion in every boundary element analysis. 

We optimise the shape of the hollow, preserving the initial square-shape of both top aperture and bottom surface. To this end, we choose $32$ CPs on the side walls of the hollow; 20 of them in the back side are shown as red points in Figure~\ref{fig:app7f_config} (right). We allow each CP to move up to a certain value from its initial position. The value is selected as $0.15$, which is less than the half of the minimum distance (i.e. $0.4$) between any two CPs at the initial configuration, so that any CP does not touch the others as well as all the observation points.

\nprounddigits{2}

\begin{figure}[H]
  \centering
  \begin{tabular}{cc}
    \begin{minipage}{.45\textwidth}
      \begin{tikzpicture}[scale=.7]
  \begin{axis}[
    xlabel={Evaluation count},
    ylabel style={align=center},
    ylabel={Objective function $\mathcalJ$},
    xmin=0, xmax=10, xtick={0,2,...,10},
    ymin=0.0, ymax=3.0, ytick={0.0,0.5,...,3.0},
    xmajorgrids=true,
    ymajorgrids=true,
    legend pos=north west,
    ]
    \addplot table [sharp plot, x=count, y=f] {data/app7f_nlopt40_CHECK_COUNT.txt};
    \addlegendentry{$-z$ incident}
  \end{axis}
\end{tikzpicture}
    \end{minipage}
    & 
      \begin{minipage}{.45\textwidth}
        \begin{tikzpicture}[scale=.7]
  \begin{axis}[
    xlabel={Evaluation count},
    ylabel style={align=center},
    ylabel={Objective function $\mathcalJ$},
    xmin=0, xmax=60, xtick={0,10,...,60},
    ymin=0.0, ymax=50.0, ytick={0.0,10.0,...,50.0},
    xmajorgrids=true,
    ymajorgrids=true,
    legend pos=north west,
    ]
    \addplot [red, mark=square*] table [sharp plot, x=count, y=f] {data/app7g_nlopt40_CHECK_COUNT.txt};
    \addlegendentry{$-x$ incident}
  \end{axis}
\end{tikzpicture}
      \end{minipage}
  \end{tabular}
  \caption{History of the value of the objective function $\mathcalJ$ in the resonator model (Section~\ref{s:resonator}).}
  \label{fig:app7fg_history}
\end{figure}
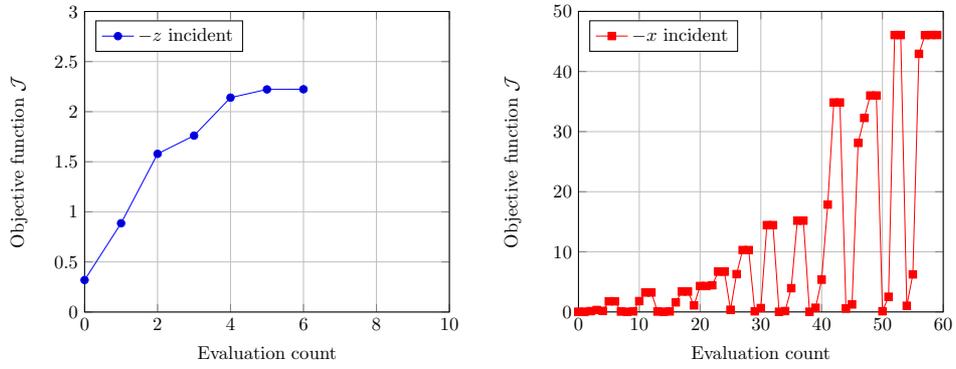

Figure~\ref{fig:app7fg_history} shows the history of $\mathcalJ$ in both incident fields. We could increase $\mathcalJ$ in both cases. The value of $\mathcalJ$ increased monotonically from \numprint{3.19436966151474544e-01} and then converged at \numprint{2.22377289810397105e+00} after seven counts in the case of the vertical ($-z$-direction) incidence, while it oscillated significantly but increase gradually from \numprint{4.12797175555163764e-02} to \numprint{4.60735738700995157e+01} after 60 counts in the case of the horizontal ($-x$-direction) incidence.

\begin{figure}[H]
  \begin{tabular}{cc}
    Vertical incidence: Initial shape & Vertical incidence: Optimised shape\\
    \includegraphics[width=.48\textwidth]{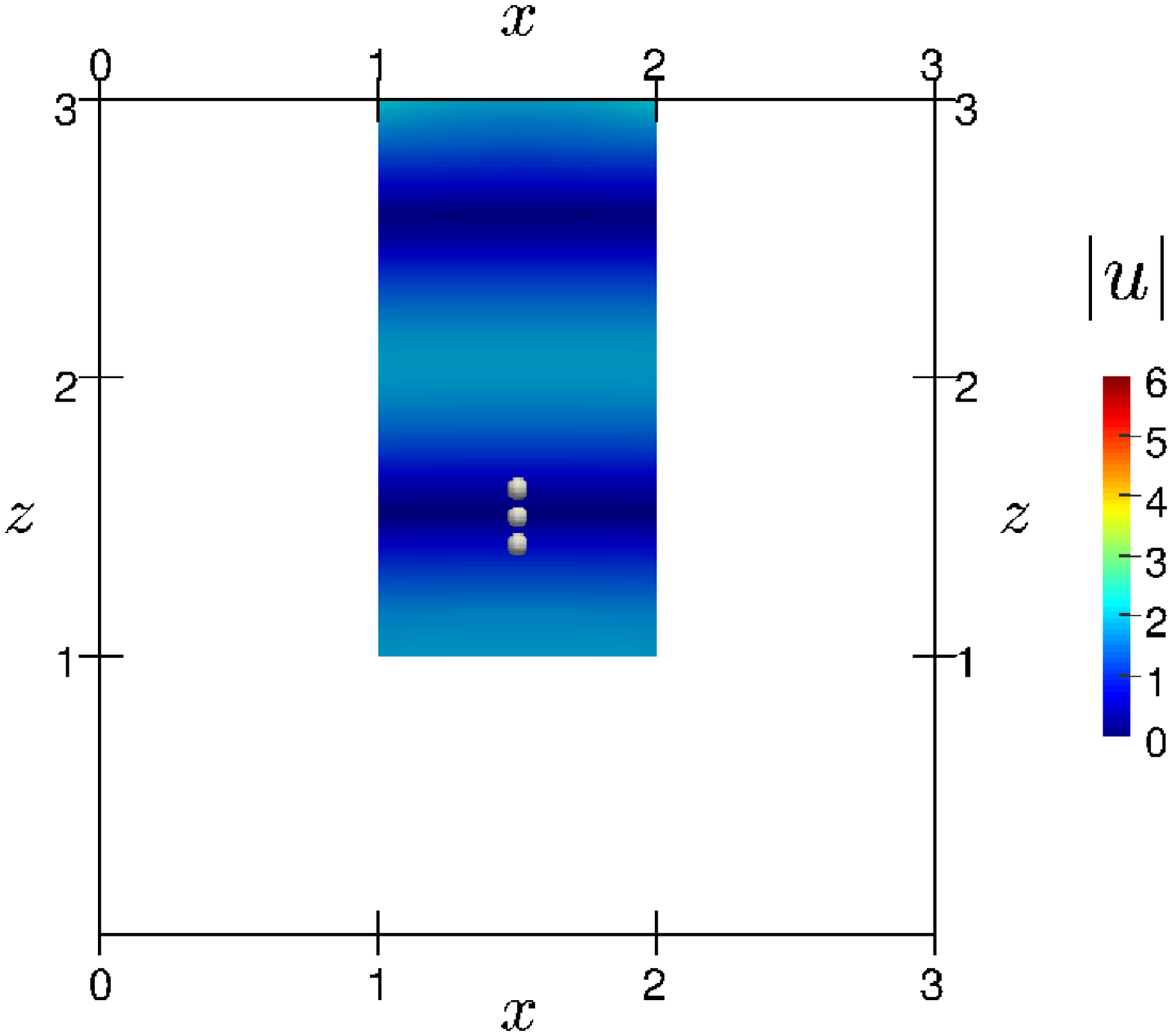}
    &\includegraphics[width=.48\textwidth]{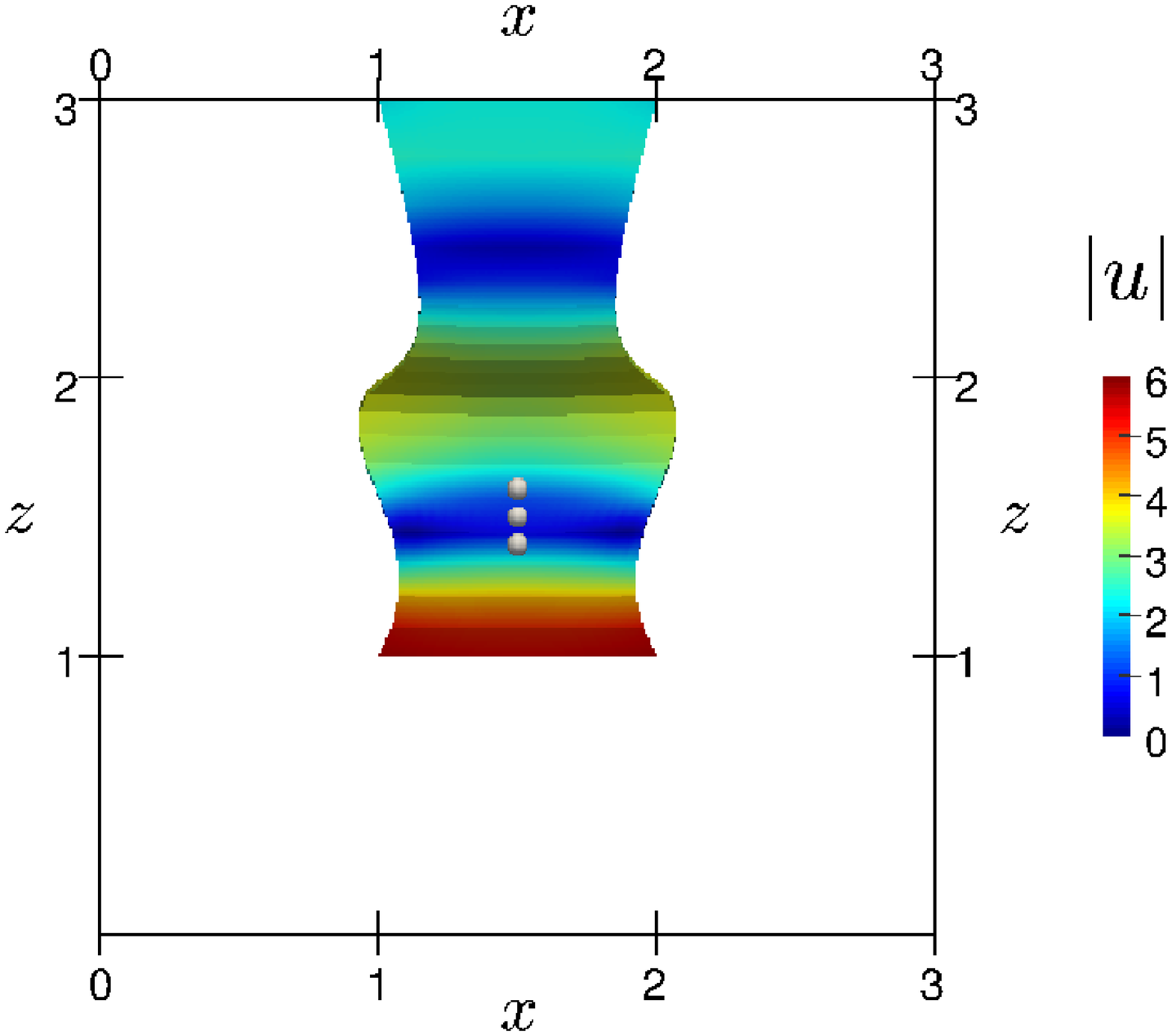}\\
    Horizontal incidence: Initial shape & Horizontal incidence: Optimised shape\\
    \includegraphics[width=.48\textwidth]{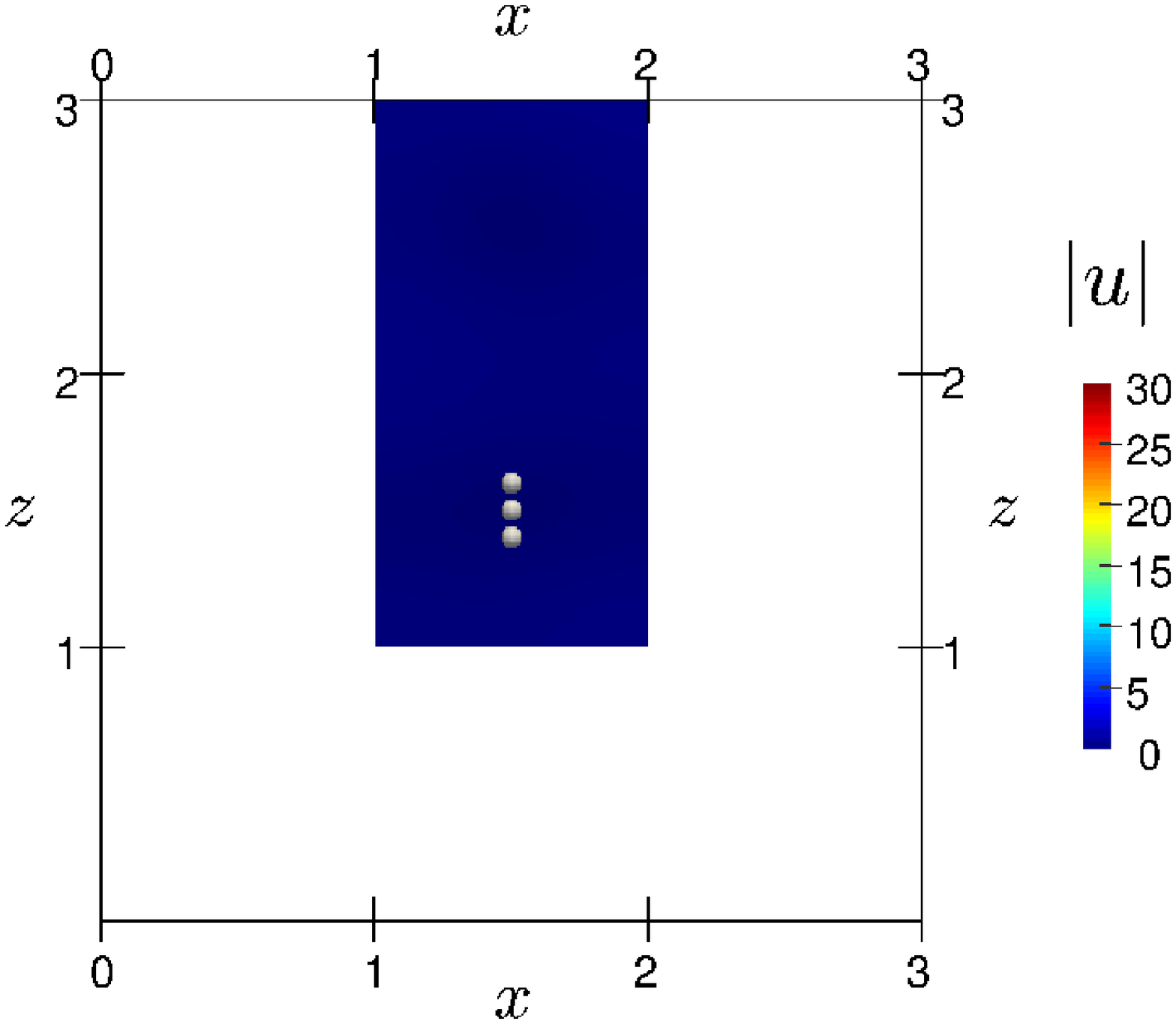}
    &\includegraphics[width=.48\textwidth]{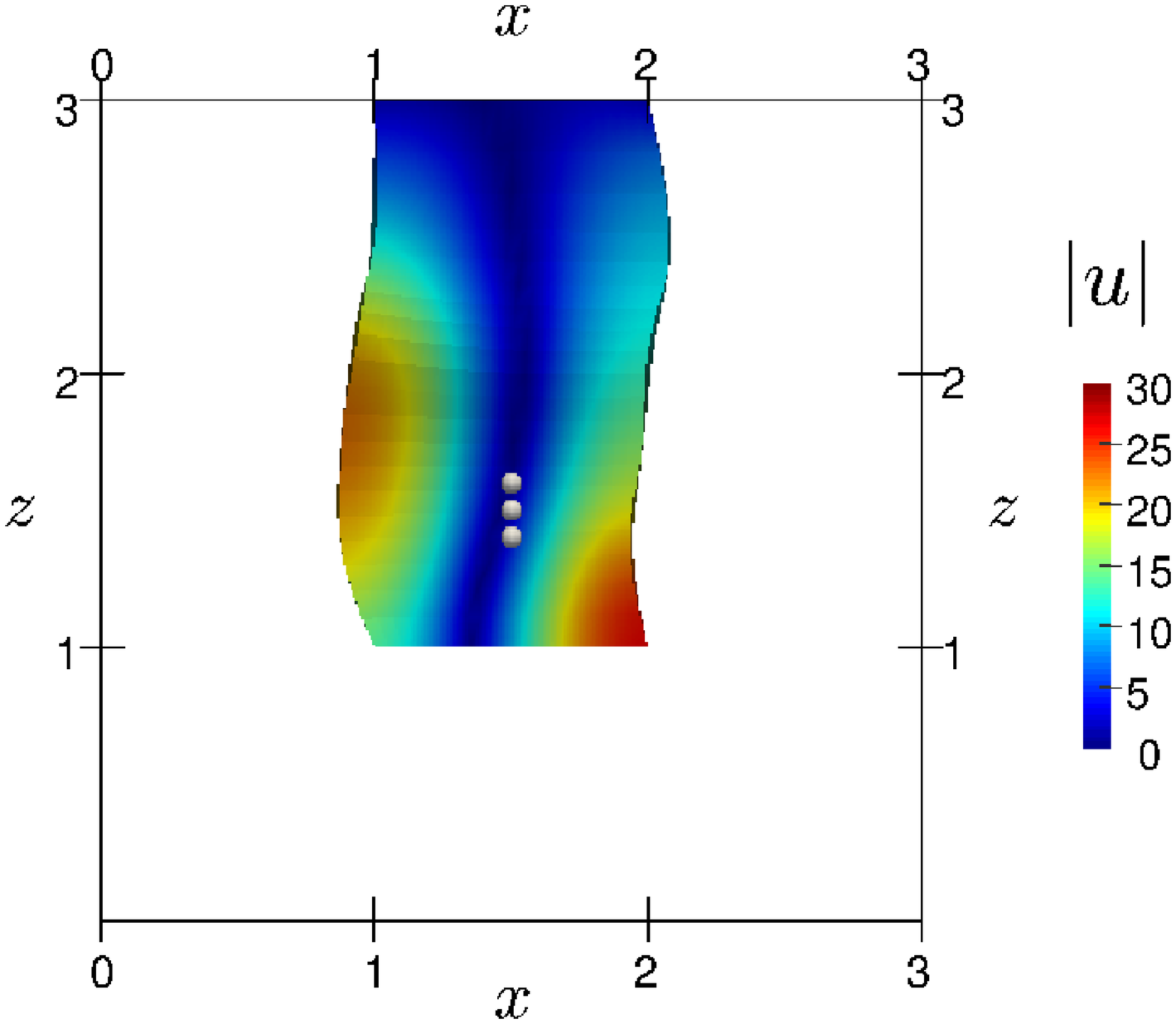}
  \end{tabular}
  \caption{Initial (left) and optimised (right) shapes of the hollow together with the distribution of $|u|$ in the case of the horizontal and horizontal incidences (Section~\ref{s:resonator}). Here, the range of $y$ is from 1.5 to 2.5, which covers the back half of the hollow. The points in grey are the observation points.}
  \label{fig:app7_u_abs_y15_25}
\end{figure}

\begin{figure}[H]
  \centering
  \includegraphics[width=.8\textwidth]{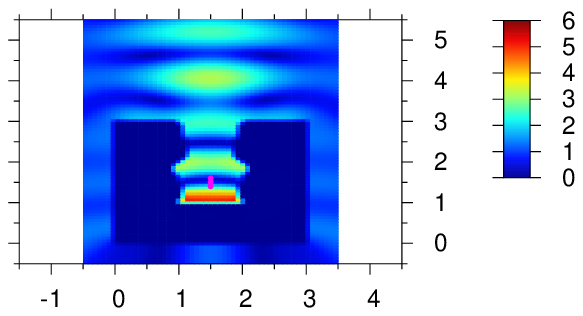}\\
  Vertical incident planewave\\
  \includegraphics[width=.8\textwidth]{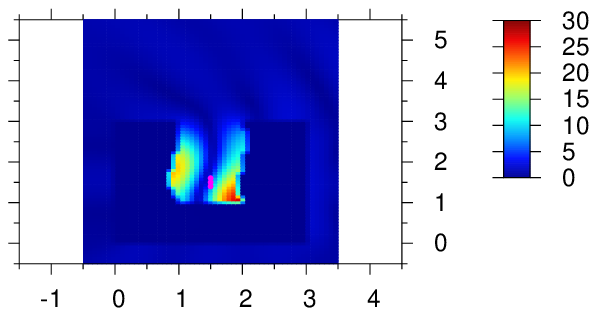}\\
  Horizontal incident planewave
  \caption{Distribution of $|u|$ on the middle cross section of $y=1.5$ in the optimised shape in the case of the vertical and horizontal incidences (Section~\ref{s:resonator}). The points coloured in magenta represents the observation points.}
  \label{fig:app7fg_uin_final}
\end{figure}

Figures~\ref{fig:app7_u_abs_y15_25} shows the initial and optimised (final) shapes of the hollow together with the distribution of $|u|$ on it. In both cases, the $|u|$ around the observation points was relatively low, but the sound pressure was actually strengthened inside the hollow after each optimisation. In addition, Figure~\ref{fig:app7fg_uin_final} shows the distribution of $|u|$ on the middle cross section of $y=1.5$.

\subsection{Example 3: Bending duct}\label{s:bending}

\begin{figure}[H]
  \centering
  \includegraphics[width=.55\textwidth]{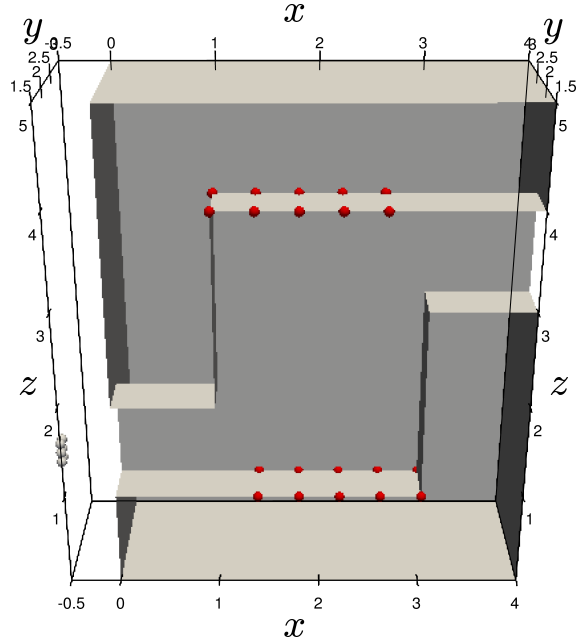}
  \caption{The back ($-y$) side of the bending-duct model at the initial configuration (Section~\ref{s:bending}). The points in grey on the plane $x=-0.5$ represent six of the nine observation point $\{\bm{z}_m\}_{m=1}^9$, which are close to the exit of the duct on the plane $x=0$. The red points represent 20 of the 30 CPs to be designed.}
  \label{fig:app17_config}
\end{figure}

As the final example, we consider a more complicated model, whose dimension is $3\times 3 \times 5$, containing a bending duct. Figure~\ref{fig:app17_config} shows the cross section of $y=1.5$ to see the inside of the model. As illustrated in the figure, we consider $3\times 3$ observation points on the plane $x=-0.5$ near the exit of the duct. We design a part of the top and bottom surfaces of the duct through 30 control points, 20 of which are drawn as red points in the figure. More precisely, we optimise the $z$ coordinates of those CPs. Here, every coordinate can be changed from its initial value up to 0.2, which takes account of the aforementioned consideration, that is, any CP never collides with the others during the optimisation. The incident planewave is given from the $+x$ side and its wavenumber is selected as 1, 2, 3 or 5 for comparison.

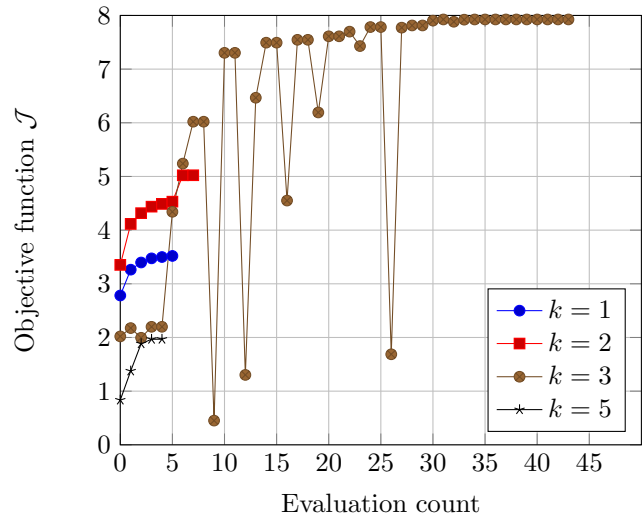
\begin{figure}[H]
  \centering
  \begin{tikzpicture}[scale=1.0]
  \begin{axis}[
    xlabel={Evaluation count},
    ylabel style={align=center},
    ylabel={Objective function $\mathcalJ$},
    xmin=0, xmax=50, xtick={0,5,...,45},
    ymin=0.0, ymax=8.0, ytick={0.0,1.0,...,8.0},
    xmajorgrids=true,
    ymajorgrids=true,
    legend pos=south east,
    ]
    \addplot table [sharp plot, x=count, y=f] {data/app17c_nlopt40_CHECK_COUNT.txt};
    \addlegendentry{$k=1$}

    \addplot table [sharp plot, x=count, y=f] {data/app17b_nlopt40_CHECK_COUNT.txt};
    \addlegendentry{$k=2$}

    \addplot table [sharp plot, x=count, y=f] {data/app17a_nlopt40_CHECK_COUNT.txt};
    \addlegendentry{$k=3$}

    \addplot table [sharp plot, x=count, y=f] {data/app17d_nlopt40_CHECK_COUNT.txt};
    \addlegendentry{$k=5$}

  \end{axis}
\end{tikzpicture}
  \caption{History of the value of the objective function $\mathcalJ$ in the bending-duct model (Section~\ref{s:bending}).}
  \label{fig:app17_history}
\end{figure}

\begin{figure}[H]
  \centering
  \begin{tabular}{ccc}
    & Initial &  Optimised\\
    $k=1$
    &\includegraphics[height=.22\textheight]{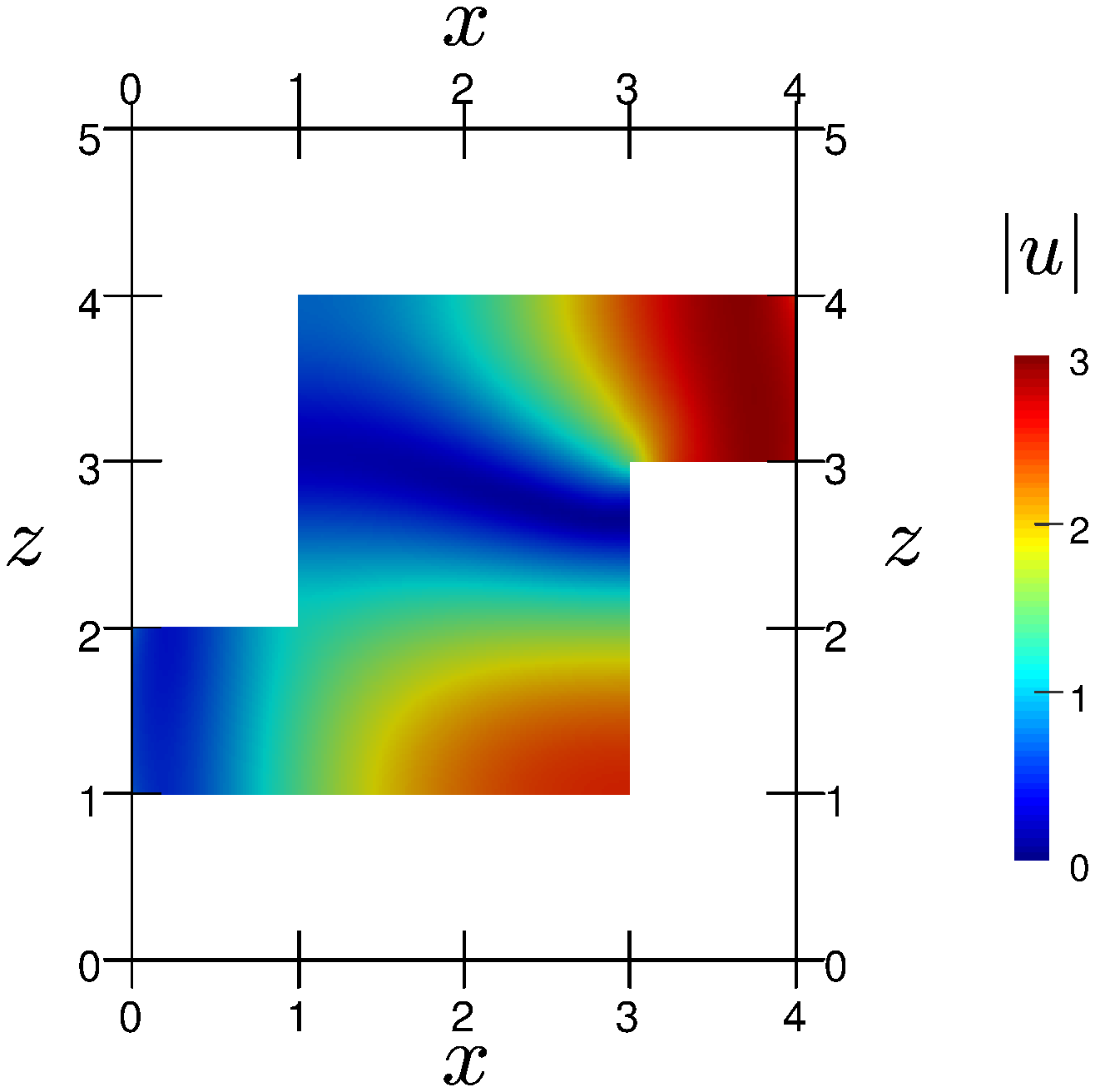}
    &\includegraphics[height=.22\textheight]{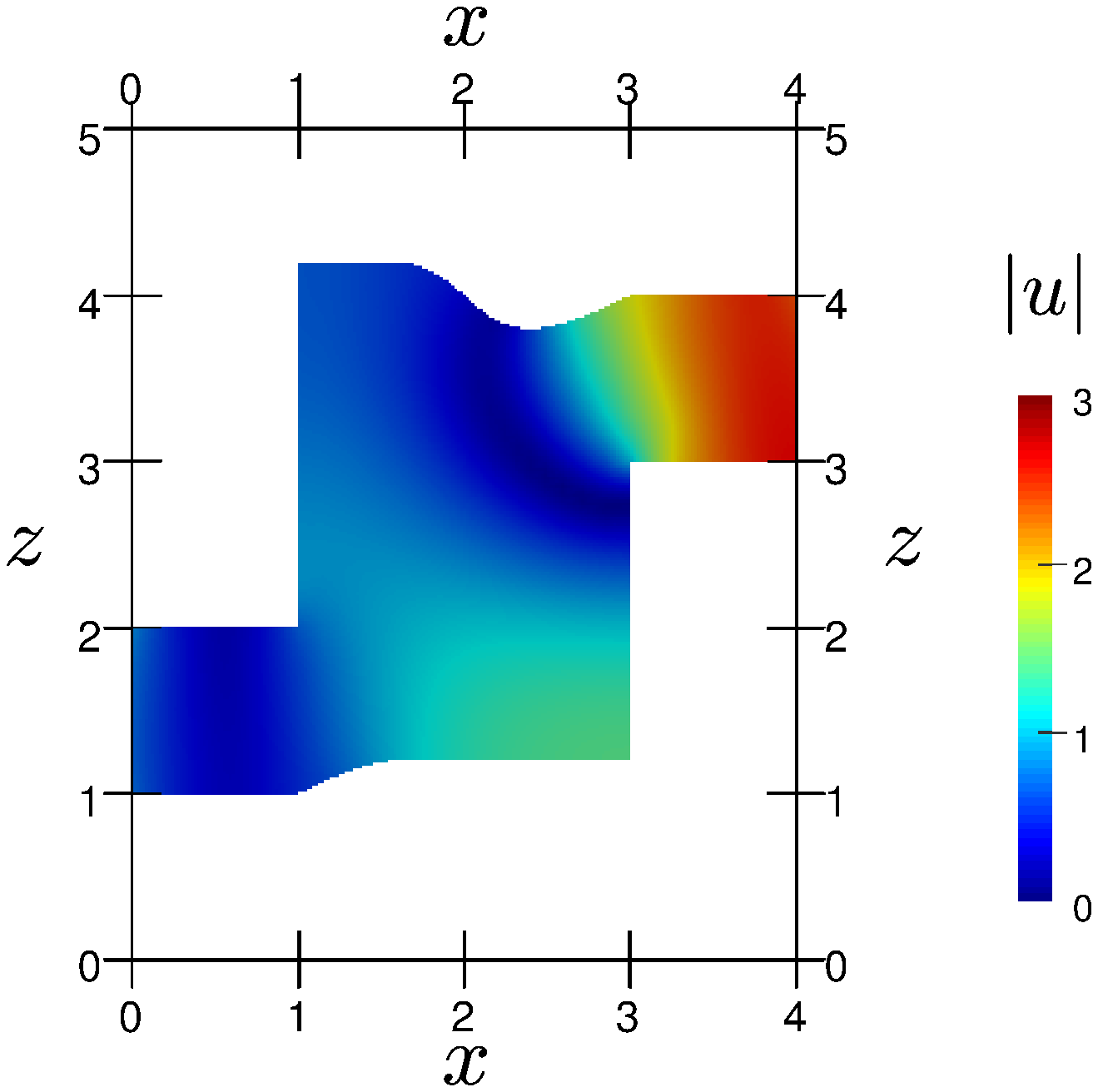}\\
    $k=2$
    &\includegraphics[height=.22\textheight]{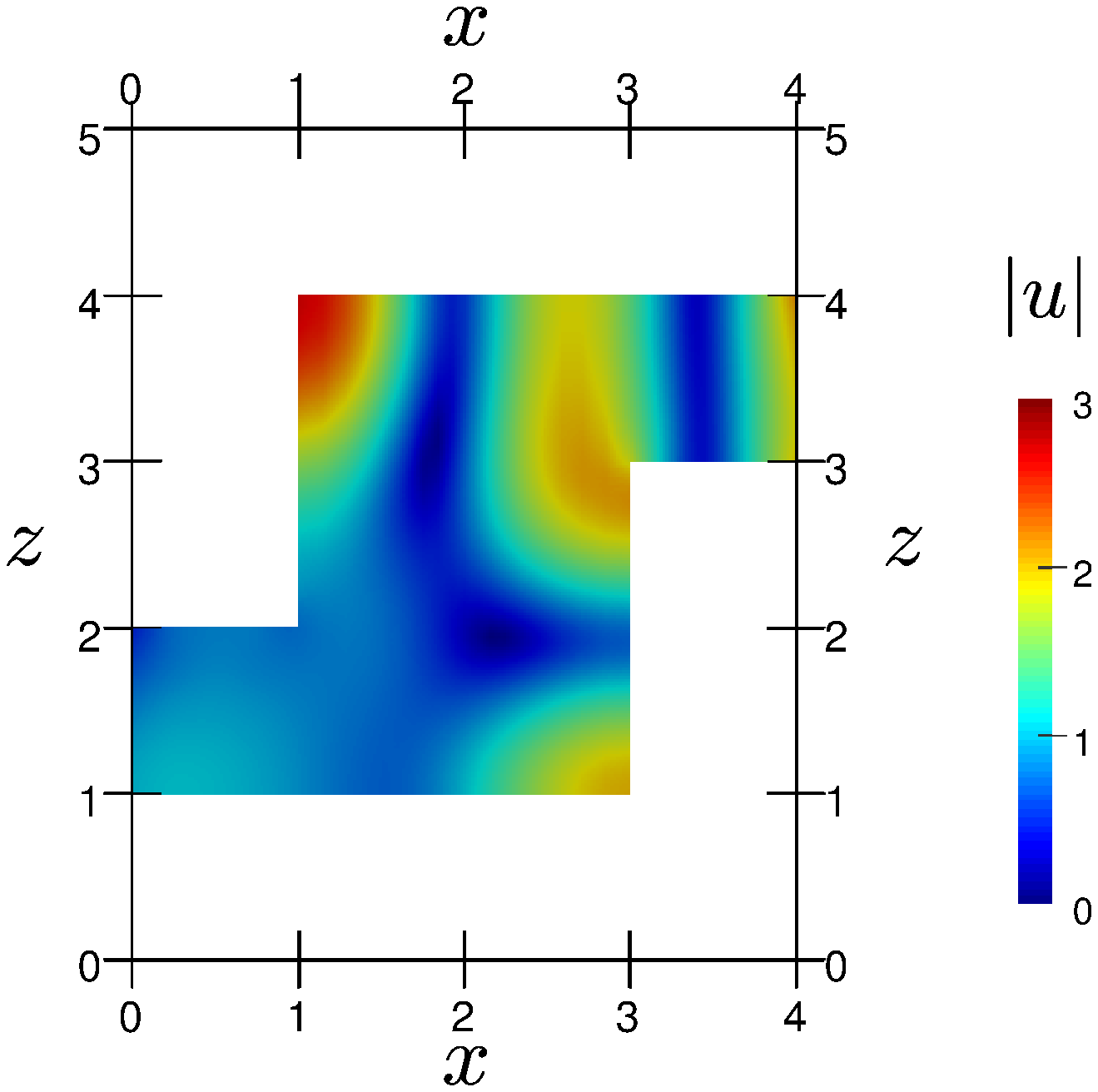}
    &\includegraphics[height=.22\textheight]{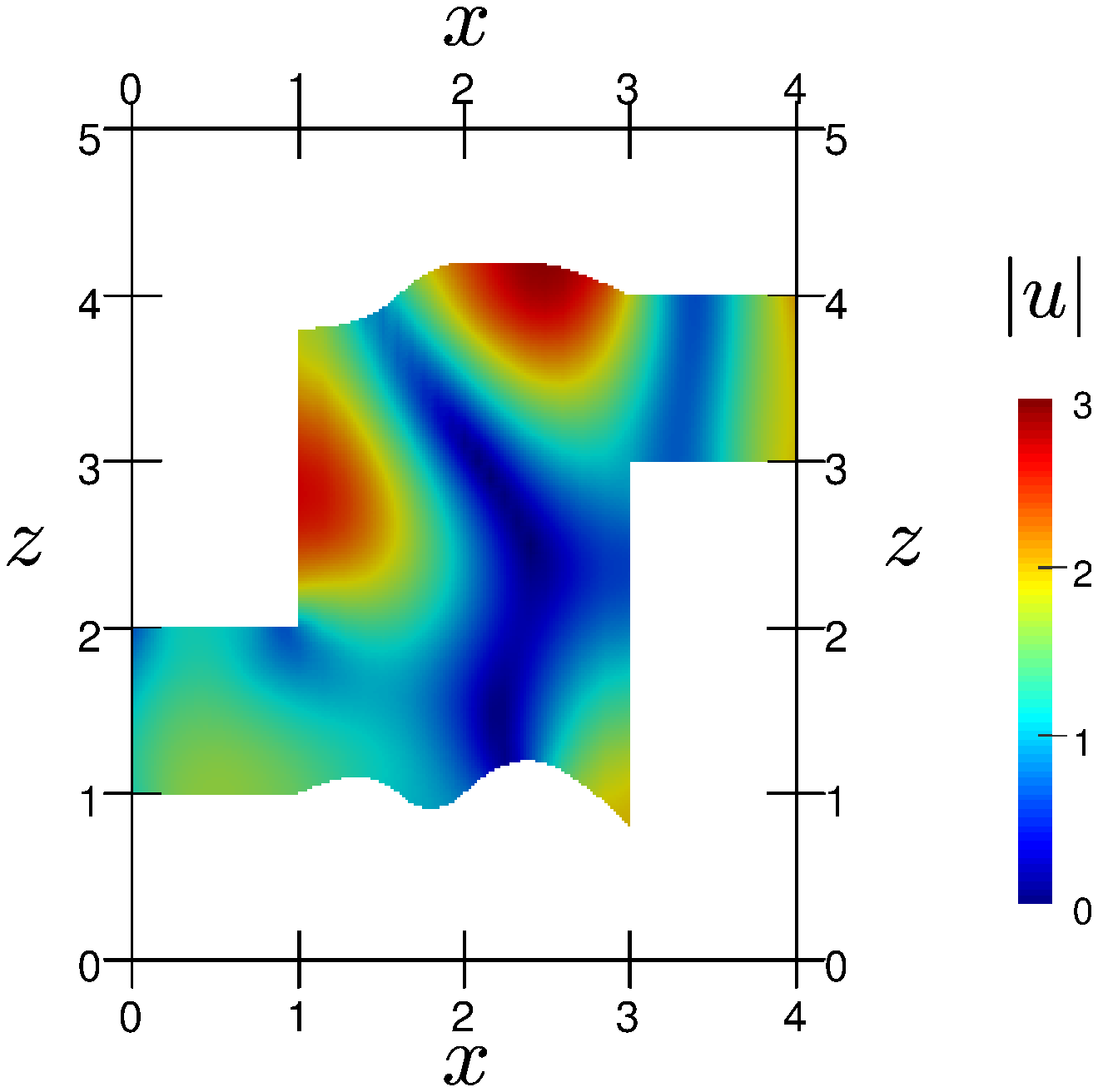}\\
    $k=3$
    &\includegraphics[height=.22\textheight]{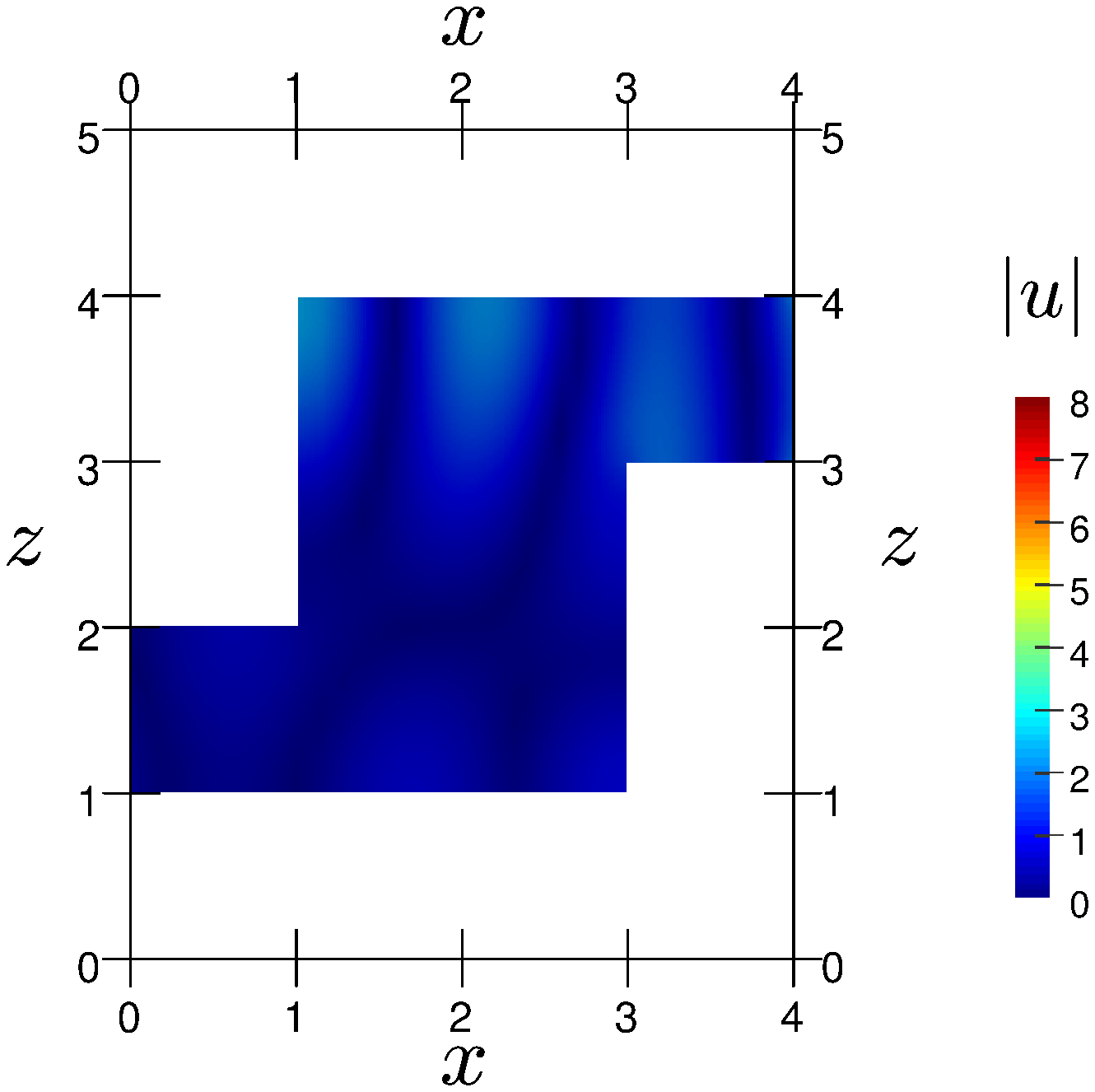}
    &\includegraphics[height=.22\textheight]{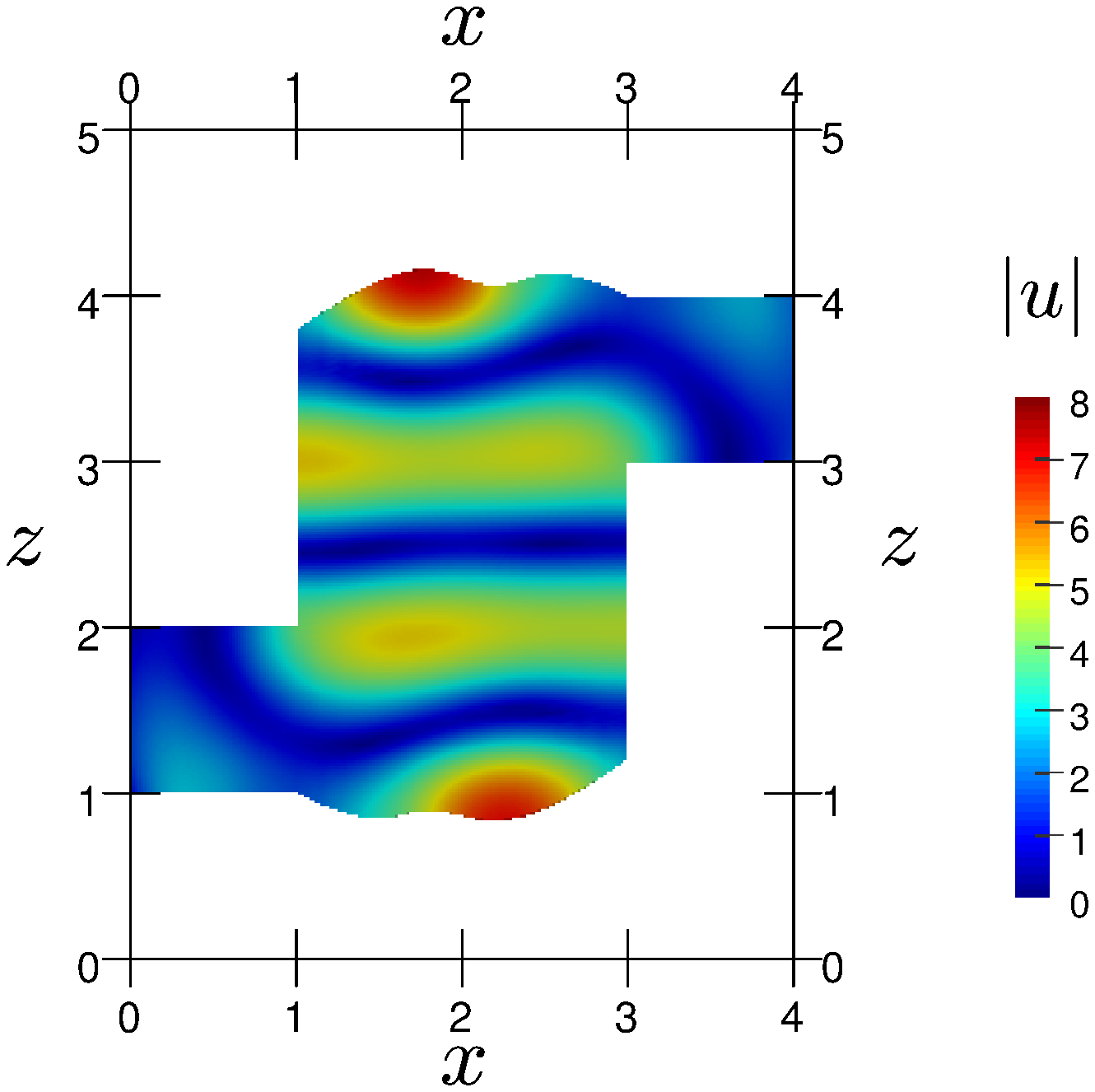}\\
    $k=5$
    &\includegraphics[height=.22\textheight]{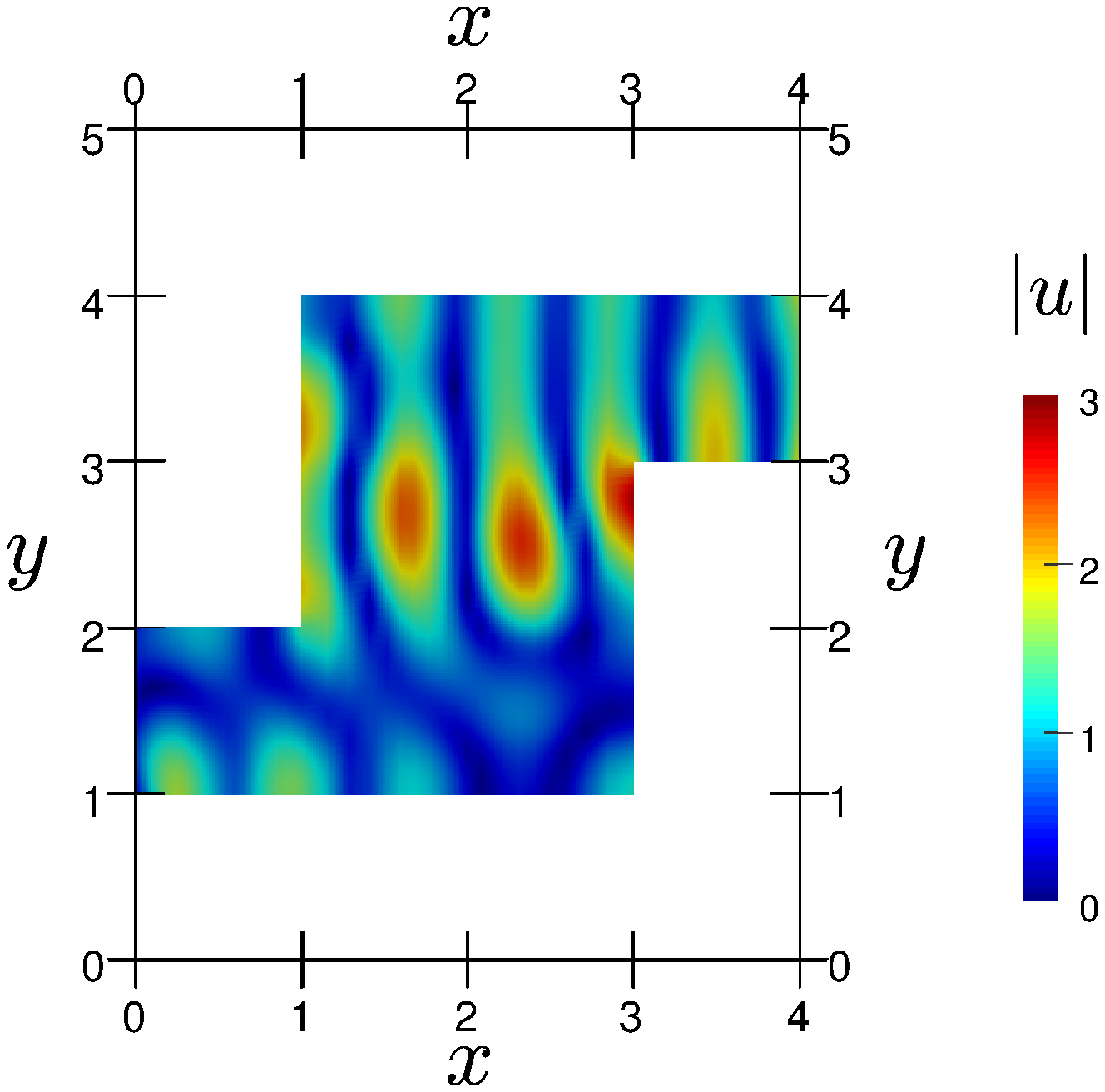}
    &\includegraphics[height=.22\textheight]{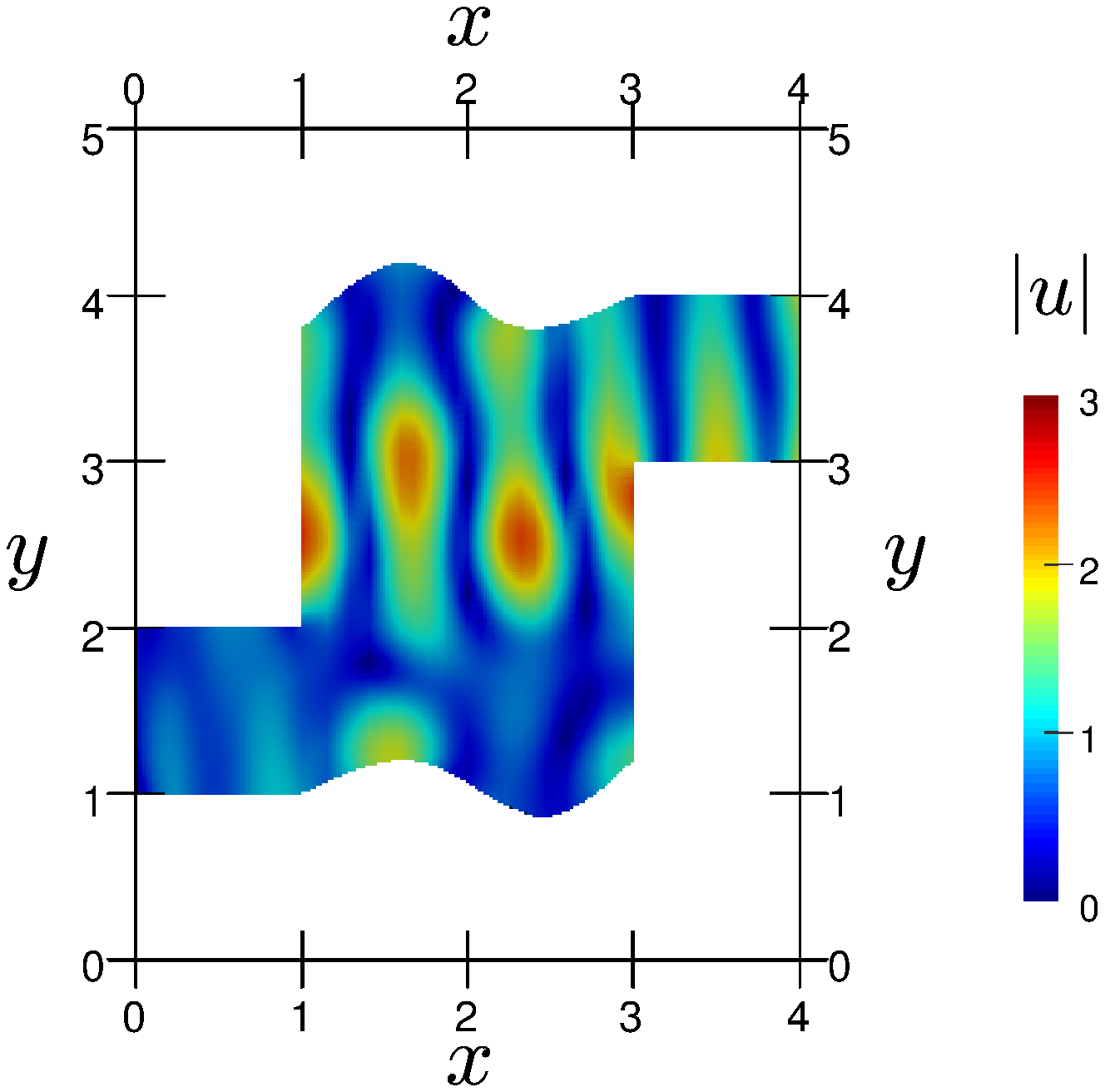}
  \end{tabular}      
  \caption{Shape and distribution of $|u|$ at the initial (left) and optimised (right) configurations of the bending-duct model (Section~\ref{s:bending}). Note that the range of $|u|$ varies according to the wavenumber $k$.}
  \label{fig:app17}
\end{figure}

As shown in Figure~\ref{fig:app17_history}, the optimisation was successfully terminated for every wavenumber $k$. The value of $\mathcalJ$ was largest at the second largest $k=3$. This can be related to the standing wave in the vertical direction excited in the duct, which can be observed in Figure~\ref{fig:app17}, but we did not pursue the reason.

\section{Conclusion}\label{s:concl}

Exploiting our previous work~\cite{takahashi2018jascome} on the development of the isogeometric boundary element method (IGBEM) for the 3D Helmholtz equation and the sensitivity analysis based on the adjoint variable method (AVM), we have newly proposed a shape optimisation system by integrating the IGBEM and AVM into nonlinear optimisation algorithms, viz. a primal-dual interior point method (IP) and the method of moving asymptotes (MMA) as well as the sequential least-squares quadratic programming (SLSQP), which are available from the open software Ipopt~\cite{ipopt} and NLopt~\cite{NLopt}, respectively. We have numerically verified the system in a (parametric) optimisation problem that has the exact solution. The system could find the optimal solutions successfully. Then, we applied the system to optimise three models that consider a reflector, resonator and bending-duct, which consist of multiple NURBS surfaces. We could maximise the objective function in every optimisation and found that the SLSQP was the best in the sense that it required less number of evaluating the objective function as well as its gradient, which is the most time-consuming part in the optimisation based on the IGBEM and AVM.

In the future, we are going to enhance our shape optimisation system so that it can directly deal with NURBS models that are generated by solid modellers such as Rhinoceras\footnote{Home page: \texttt{https://www.rhino3d.com/}.} and SMlib\footnote{Home page: \texttt{https://smlib.com/smlib/}.}. This task is technically evident but practically important to give initial shapes comprising of truly curved surfaces. In addition, we are planning to develop a similar shape optimisation software for electromagnetism in 3D from the present one, considering the application to metamaterials and plasmonics.

\section*{Acknowledgements}

This research was supported by JSPS KAKENHI (Grant number: 18K11335).

\section*{\refname}
\iffalse
\bibliographystyle{elsarticle/elsarticle-num}
\bibliography{ref,eat}

\begin{thebibliography}{10}
\expandafter\ifx\csname url\endcsname\relax
  \def\url#1{\texttt{#1}}\fi
\expandafter\ifx\csname urlprefix\endcsname\relax\def\urlprefix{URL }\fi
\expandafter\ifx\csname href\endcsname\relax
  \def\href#1#2{#2} \def\path#1{#1}\fi

\bibitem{Ebbesen_1998}
T.~W. Ebbesen, H.~J. Lezec, H.~F. Ghaemi, T.~Thio, P.~A. Wolff,
  \href{http://dx.doi.org/10.1038/35570}{Extraordinary optical transmission
  through sub-wavelength hole arrays}, Nature 391~(6668) (1998) 667--669.
\newblock \href {http://dx.doi.org/10.1038/35570} {\path{doi:10.1038/35570}}.
\newline\urlprefix\url{http://dx.doi.org/10.1038/35570}

\bibitem{Lezec_2002}
H.~J. Lezec, A.~Degiron, E.~Devaux, R.~A. Linke, L.~Martin-Moreno, F.~J.
  Garcia-Vidal, T.~W. Ebbesen,
  \href{http://www.sciencemag.org/content/297/5582/820.abstract}{Beaming light
  from a subwavelength aperture}, Science 297~(5582) (2002) 820--822.
\newblock \href {http://dx.doi.org/10.1126/science.1071895}
  {\path{doi:10.1126/science.1071895}}.
\newline\urlprefix\url{http://www.sciencemag.org/content/297/5582/820.abstract}

\bibitem{Christensen_2007}
J.~Christensen, A.~I. Fernandez-Dominguez, F.~De~Leon-Perez, L.~Martin-Moreno,
  F.~J. Garcia-Vidal, Collimation of sound assisted by acoustic surface waves,
  Nature Physics 3 (2007) 851--852.

\bibitem{Zhou_2010}
Y.~Zhou, M.-H. Lu, L.~Feng, X.~Ni, Y.-F. Chen, Y.-Y. Zhu, S.-N. Zhu, N.-B.
  Ming, \href{http://link.aps.org/doi/10.1103/PhysRevLett.104.164301}{Acoustic
  surface evanescent wave and its dominant contribution to extraordinary
  acoustic transmission and collimation of sound}, Phys. Rev. Lett. 104 (2010)
  164301.
\newblock \href {http://dx.doi.org/10.1103/PhysRevLett.104.164301}
  {\path{doi:10.1103/PhysRevLett.104.164301}}.
\newline\urlprefix\url{http://link.aps.org/doi/10.1103/PhysRevLett.104.164301}

\bibitem{takahashi2014}
T.~Takahashi, K.~Kuriyama, T.~Matsumoto,
  \href{https://doi.org/10.1063/1.4753801}{Beaming of inplane elastic waves
  through a subwavelength channel with periodic corrugations}, Applied Physics
  Letters 101~(12) (2012) 124101.
\newblock \href {http://dx.doi.org/10.1063/1.4753801}
  {\path{doi:10.1063/1.4753801}}.
\newline\urlprefix\url{https://doi.org/10.1063/1.4753801}

\bibitem{liu2011metamaterials}
Y.~Liu, X.~Zhang, \href{http://dx.doi.org/10.1039/C0CS00184H}{Metamaterials: a
  new frontier of science and technology}, Chem. Soc. Rev. 40 (2011)
  2494--2507.
\newblock \href {http://dx.doi.org/10.1039/C0CS00184H}
  {\path{doi:10.1039/C0CS00184H}}.
\newline\urlprefix\url{http://dx.doi.org/10.1039/C0CS00184H}

\bibitem{wang2020tunable}
W.~Yan-Feng, W.~Yi-Ze, W.~Bin, C.~Weiqiu, W.~Yue-Sheng,
  \href{https://doi.org/10.1115/1.4046222}{{Tunable and Active Phononic
  Crystals and Metamaterials}}, Applied Mechanics Reviews 72~(4), 040801.
\newblock \href {http://dx.doi.org/10.1115/1.4046222}
  {\path{doi:10.1115/1.4046222}}.
\newline\urlprefix\url{https://doi.org/10.1115/1.4046222}

\bibitem{soares1984}
S.~C., R.~H\'{e}lder, F.~Luis, H.~Edward, Optimization of the geometry of
  shafts using boundary elements, Journal of Mechanisms Transmissions and
  Automation in Design 106 (1984) 199--202.
\newblock \href {http://dx.doi.org/10.1115/1.3258579}
  {\path{doi:10.1115/1.3258579}}.

\bibitem{hughes2005}
T.~Hughes, J.~Cottrell, Y.~Bazilevs,
  \href{http://www.sciencedirect.com/science/article/pii/S0045782504005171}{Isogeometric
  analysis: {CAD}, finite elements, {NURBS}, exact geometry and mesh
  refinement}, Computer Methods in Applied Mechanics and Engineering
  194~(39--41) (2005) 4135--4195.
\newblock \href {http://dx.doi.org/10.1016/j.cma.2004.10.008}
  {\path{doi:10.1016/j.cma.2004.10.008}}.
\newline\urlprefix\url{http://www.sciencedirect.com/science/article/pii/S0045782504005171}

\bibitem{cottrell2009}
J.~Cottrell, T.~Hughes, Y.~Bazilevs, Isogeometric Analysis: Toward integration
  of CAD and FEA, 2009.
\newblock \href {http://dx.doi.org/10.1002/9780470749081.ch7}
  {\path{doi:10.1002/9780470749081.ch7}}.

\bibitem{yoon2015}
M.~Yoon, M.-J. Choi, S.~Cho,
  \href{http://www.sciencedirect.com/science/article/pii/S001793101500616X}{Isogeometric
  configuration design optimization of heat conduction problems using boundary
  integral equation}, International Journal of Heat and Mass Transfer 89 (2015)
  937--949.
\newblock \href
  {http://dx.doi.org/http://dx.doi.org/10.1016/j.ijheatmasstransfer.2015.05.112}
  {\path{doi:http://dx.doi.org/10.1016/j.ijheatmasstransfer.2015.05.112}}.
\newline\urlprefix\url{http://www.sciencedirect.com/science/article/pii/S001793101500616X}

\bibitem{kostas2015}
K.~Kostas, A.~Ginnis, C.~Politis, P.~Kaklis,
  \href{http://www.sciencedirect.com/science/article/pii/S0045782514004009}{Ship-hull
  shape optimization with a {T-spline} based {BEM}-isogeometric solver},
  Computer Methods in Applied Mechanics and Engineering 284 (2015) 611--622,
  isogeometric Analysis Special Issue.
\newblock \href {http://dx.doi.org/http://dx.doi.org/10.1016/j.cma.2014.10.030}
  {\path{doi:http://dx.doi.org/10.1016/j.cma.2014.10.030}}.
\newline\urlprefix\url{http://www.sciencedirect.com/science/article/pii/S0045782514004009}

\bibitem{gillebaart2016}
E.~Gillebaart, R.~D. Breuker,
  \href{http://www.sciencedirect.com/science/article/pii/S0045782516300962}{Low-fidelity
  {2D} isogeometric aeroelastic analysis and optimization method with
  application to a morphing airfoil}, Computer Methods in Applied Mechanics and
  Engineering 305 (2016) 512--536.
\newblock \href {http://dx.doi.org/http://dx.doi.org/10.1016/j.cma.2016.03.014}
  {\path{doi:http://dx.doi.org/10.1016/j.cma.2016.03.014}}.
\newline\urlprefix\url{http://www.sciencedirect.com/science/article/pii/S0045782516300962}

\bibitem{kostas2017}
K.~Kostas, A.~Ginnis, C.~Politis, P.~Kaklis,
  \href{http://www.sciencedirect.com/science/article/pii/S0010448516300653}{Shape-optimization
  of {2D} hydrofoils using an isogeometric {BEM} solver}, Computer-Aided Design
  82 (2017) 79--87, isogeometric Design and Analysis.
\newblock \href {http://dx.doi.org/https://doi.org/10.1016/j.cad.2016.07.002}
  {\path{doi:https://doi.org/10.1016/j.cad.2016.07.002}}.
\newline\urlprefix\url{http://www.sciencedirect.com/science/article/pii/S0010448516300653}

\bibitem{kostas2018}
K.~Kostas, M.~Fyrillas, C.~Politis, A.~Ginnis, P.~Kaklis,
  \href{https://www.sciencedirect.com/science/article/pii/S0045782518303177}{Shape
  optimization of conductive-media interfaces using an {IGA-BEM} solver},
  Computer Methods in Applied Mechanics and Engineering 340 (2018) 600--614.
\newblock \href {http://dx.doi.org/https://doi.org/10.1016/j.cma.2018.06.019}
  {\path{doi:https://doi.org/10.1016/j.cma.2018.06.019}}.
\newline\urlprefix\url{https://www.sciencedirect.com/science/article/pii/S0045782518303177}

\bibitem{li2011}
K.~Li, X.~Qian,
  \href{http://www.sciencedirect.com/science/article/pii/S0010448511002302}{Isogeometric
  analysis and shape optimization via boundary integral}, Computer-Aided Design
  43~(11) (2011) 1427--1437.
\newblock \href {http://dx.doi.org/http://dx.doi.org/10.1016/j.cad.2011.08.031}
  {\path{doi:http://dx.doi.org/10.1016/j.cad.2011.08.031}}.
\newline\urlprefix\url{http://www.sciencedirect.com/science/article/pii/S0010448511002302}

\bibitem{lian2016}
H.~Lian, P.~Kerfriden, S.~P.~A. Bordas,
  \href{http://dx.doi.org/10.1002/nme.5149}{Implementation of regularized
  isogeometric boundary element methods for gradient-based shape optimization
  in two-dimensional linear elasticity}, International Journal for Numerical
  Methods in Engineering 106~(12) (2016) 972--1017, nme.5149.
\newblock \href {http://dx.doi.org/10.1002/nme.5149}
  {\path{doi:10.1002/nme.5149}}.
\newline\urlprefix\url{http://dx.doi.org/10.1002/nme.5149}

\bibitem{lian2017}
H.~Lian, P.~Kerfriden, S.~Bordas,
  \href{http://www.sciencedirect.com/science/article/pii/S0045782516315365}{Shape
  optimization directly from {CAD}: An isogeometric boundary element approach
  using {T}-splines}, Computer Methods in Applied Mechanics and Engineering 317
  (2017) 1--41.
\newblock \href {http://dx.doi.org/https://doi.org/10.1016/j.cma.2016.11.012}
  {\path{doi:https://doi.org/10.1016/j.cma.2016.11.012}}.
\newline\urlprefix\url{http://www.sciencedirect.com/science/article/pii/S0045782516315365}

\bibitem{sun2018}
S.~Sun, T.~Yu, T.~Nguyen, E.~Atroshchenko, T.~Bui,
  \href{https://www.sciencedirect.com/science/article/pii/S0955799717303053}{Structural
  shape optimization by igabem and particle swarm optimization algorithm},
  Engineering Analysis with Boundary Elements 88 (2018) 26--40.
\newblock \href
  {http://dx.doi.org/https://doi.org/10.1016/j.enganabound.2017.12.007}
  {\path{doi:https://doi.org/10.1016/j.enganabound.2017.12.007}}.
\newline\urlprefix\url{https://www.sciencedirect.com/science/article/pii/S0955799717303053}

\bibitem{li2019}
S.~Li, J.~Trevelyan, Z.~Wu, H.~Lian, D.~Wang, W.~Zhang,
  \href{https://www.sciencedirect.com/science/article/pii/S0045782519300945}{An
  adaptive {SVD}--{Krylov} reduced order model for surrogate based structural
  shape optimization through isogeometric boundary element method}, Computer
  Methods in Applied Mechanics and Engineering 349 (2019) 312--338.
\newblock \href {http://dx.doi.org/https://doi.org/10.1016/j.cma.2019.02.023}
  {\path{doi:https://doi.org/10.1016/j.cma.2019.02.023}}.
\newline\urlprefix\url{https://www.sciencedirect.com/science/article/pii/S0045782519300945}

\bibitem{sun2020}
D.~Sun, C.~Dong,
  \href{https://www.sciencedirect.com/science/article/pii/S0045782520304643}{Shape
  optimization of heterogeneous materials based on isogeometric boundary
  element method}, Computer Methods in Applied Mechanics and Engineering 370
  (2020) 113279.
\newblock \href {http://dx.doi.org/https://doi.org/10.1016/j.cma.2020.113279}
  {\path{doi:https://doi.org/10.1016/j.cma.2020.113279}}.
\newline\urlprefix\url{https://www.sciencedirect.com/science/article/pii/S0045782520304643}

\bibitem{yoon2020}
M.~Yoon, J.~Lee, B.~Koo,
  \href{https://www.sciencedirect.com/science/article/pii/S0965997820301952}{Shape
  design optimization of thermoelasticity problems using isogeometric boundary
  element method}, Advances in Engineering Software 149 (2020) 102871.
\newblock \href
  {http://dx.doi.org/https://doi.org/10.1016/j.advengsoft.2020.102871}
  {\path{doi:https://doi.org/10.1016/j.advengsoft.2020.102871}}.
\newline\urlprefix\url{https://www.sciencedirect.com/science/article/pii/S0965997820301952}

\bibitem{liu2017}
C.~Liu, L.~Chen, W.~Zhao, H.~Chen,
  \href{https://www.sciencedirect.com/science/article/pii/S0955799717300930}{Shape
  optimization of sound barrier using an isogeometric fast multipole boundary
  element method in two dimensions}, Engineering Analysis with Boundary
  Elements 85 (2017) 142--157.
\newblock \href
  {http://dx.doi.org/https://doi.org/10.1016/j.enganabound.2017.09.009}
  {\path{doi:https://doi.org/10.1016/j.enganabound.2017.09.009}}.
\newline\urlprefix\url{https://www.sciencedirect.com/science/article/pii/S0955799717300930}

\bibitem{takahashi2019ewco}
T.~Takahashi, T.~Yamamoto, Y.~Shimba, H.~Isakari, T.~Matsumoto,
  \href{https://doi.org/10.1007/s00366-018-0606-6}{A framework of shape
  optimisation based on the isogeometric boundary element method toward
  designing thin-silicon photovoltaic devices}, Eng. with Comput. 35~(2) (2019)
  423--449.
\newblock \href {http://dx.doi.org/10.1007/s00366-018-0606-6}
  {\path{doi:10.1007/s00366-018-0606-6}}.
\newline\urlprefix\url{https://doi.org/10.1007/s00366-018-0606-6}

\bibitem{ummidivarapu2020}
V.~K. Ummidivarapu, H.~K. Voruganti, T.~Khajah, S.~P.~A. Bordas,
  \href{https://www.sciencedirect.com/science/article/pii/S0167839620300686}{Isogeometric
  shape optimization of an acoustic horn using the teaching-learning-based
  optimization ({TLBO}) algorithm}, Computer Aided Geometric Design 80 (2020)
  101881.
\newblock \href {http://dx.doi.org/https://doi.org/10.1016/j.cagd.2020.101881}
  {\path{doi:https://doi.org/10.1016/j.cagd.2020.101881}}.
\newline\urlprefix\url{https://www.sciencedirect.com/science/article/pii/S0167839620300686}

\bibitem{shaaban2020}
A.~M. Shaaban, C.~Anitescu, E.~Atroshchenko, T.~Rabczuk,
  \href{https://www.sciencedirect.com/science/article/pii/S0022460X20304296}{Isogeometric
  boundary element analysis and shape optimization by {PSO} for {3D}
  axi-symmetric high frequency {Helmholtz} acoustic problems}, Journal of Sound
  and Vibration 486 (2020) 115598.
\newblock \href {http://dx.doi.org/https://doi.org/10.1016/j.jsv.2020.115598}
  {\path{doi:https://doi.org/10.1016/j.jsv.2020.115598}}.
\newline\urlprefix\url{https://www.sciencedirect.com/science/article/pii/S0022460X20304296}

\bibitem{shaaban2020b}
A.~{Mostafa Shaaban}, C.~Anitescu, E.~Atroshchenko, T.~Rabczuk,
  \href{https://www.sciencedirect.com/science/article/pii/S0955799719306733}{Shape
  optimization by conventional and extended isogeometric boundary element
  method with {PSO} for two-dimensional {Helmholtz} acoustic problems},
  Engineering Analysis with Boundary Elements 113 (2020) 156--169.
\newblock \href
  {http://dx.doi.org/https://doi.org/10.1016/j.enganabound.2019.12.012}
  {\path{doi:https://doi.org/10.1016/j.enganabound.2019.12.012}}.
\newline\urlprefix\url{https://www.sciencedirect.com/science/article/pii/S0955799719306733}

\bibitem{wang2020}
J.~Wang, C.~Zheng, L.~Chen, H.~Chen,
  \href{https://doi.org/10.1142/S2591728520500152}{Acoustic shape optimization
  based on isogeometric wideband fast multipole boundary element method with
  adjoint variable method}, Journal of Theoretical and Computational Acoustics
  28~(02) (2020) 2050015.
\newblock \href {http://dx.doi.org/10.1142/S2591728520500152}
  {\path{doi:10.1142/S2591728520500152}}.
\newline\urlprefix\url{https://doi.org/10.1142/S2591728520500152}

\bibitem{chen2019}
L.~Chen, H.~Lian, Z.~Liu, H.~Chen, E.~Atroshchenko, S.~Bordas,
  \href{https://www.sciencedirect.com/science/article/pii/S004578251930355X}{Structural
  shape optimization of three dimensional acoustic problems with isogeometric
  boundary element methods}, Computer Methods in Applied Mechanics and
  Engineering 355 (2019) 926--951.
\newblock \href {http://dx.doi.org/https://doi.org/10.1016/j.cma.2019.06.012}
  {\path{doi:https://doi.org/10.1016/j.cma.2019.06.012}}.
\newline\urlprefix\url{https://www.sciencedirect.com/science/article/pii/S004578251930355X}

\bibitem{takahashi2018jascome}
T.~Takhashi, T.~Hirai, D.~Sato, H.~Isakari, T.~Matsumoto,
  \href{http://www.matsumoto.nuem.nagoya-u.ac.jp/jascome/denshi-journal/18/JA1815.pdf}{A
  shape sensitivity analysis based on an isogeometric boundary element method
  for 3d acoustic problems}, Transactions of the Japan Society for
  Computational Methods in Engineering 18 (2018) 35--40, written in Japanese.
\newline\urlprefix\url{http://www.matsumoto.nuem.nagoya-u.ac.jp/jascome/denshi-journal/18/JA1815.pdf}

\bibitem{ipopt}
Y.~Kawajir, C.~Laird, S.~Vigerske, A.~W{\"a}chter, A tutorial for downloading,
  installing, and using {Ipopt} (revision: 2538) (February 2015).

\bibitem{NLopt}
S.~G. Johnson, \href{http://github.com/stevengj/nlopt}{{The NLopt
  nonlinear-optimization package}}.
\newline\urlprefix\url{http://github.com/stevengj/nlopt}

\bibitem{feijoo2003}
G.~R. Feij{\'{o}}o, A.~A. Oberai, P.~M. Pinsky,
  \href{https://doi.org/10.1088/0266-5611/20/1/012}{An application of shape
  optimization in the solution of inverse acoustic scattering problems},
  Inverse Problems 20~(1) (2003) 199--228.
\newblock \href {http://dx.doi.org/10.1088/0266-5611/20/1/012}
  {\path{doi:10.1088/0266-5611/20/1/012}}.
\newline\urlprefix\url{https://doi.org/10.1088/0266-5611/20/1/012}

\bibitem{ipopt_wiki}
\href{https://en.wikipedia.org/wiki/IPOPT}{Ipopt}.
\newline\urlprefix\url{https://en.wikipedia.org/wiki/IPOPT}

\bibitem{svanberg2002class}
K.~Svanberg, A class of globally convergent optimization methods based on
  conservative convex separable approximations, SIAM Journal on Optimization 12
  (2002) 555--573.

\bibitem{kraft1994}
D.~Kraft, \href{https://doi.org/10.1145/192115.192124}{Algorithm 733: {TOMP}
  --- {Fortran Modules for Optimal Control Calculations}}, ACM Trans. Math.
  Softw. 20~(3) (1994) 262--281.
\newblock \href {http://dx.doi.org/10.1145/192115.192124}
  {\path{doi:10.1145/192115.192124}}.
\newline\urlprefix\url{https://doi.org/10.1145/192115.192124}

\bibitem{cobb1988}
J.~E. Cobb,
  \href{https://collections.lib.utah.edu/dl_files/4e/77/4e7746dd53c79f8557272b92b47d2d407da4931a.pdf}{Tiling
  the sphere with rational {Bezier} patches}, Tech. Rep. Report TR UUCS-88-009,
  University of Utah, USA (1988).
\newline\urlprefix\url{https://collections.lib.utah.edu/dl_files/4e/77/4e7746dd53c79f8557272b92b47d2d407da4931a.pdf}

\bibitem{bowman1987}
J.~J. Bowman, T.~B.~A. Senior, P.~L.~E. Uslenghi,
  \href{http://adsabs.harvard.edu/cgi-bin/nph-bib\_query?bibcode=1987hpc..book.....B}{{Electromagnetic
  and Acoustic Scattering by Simple Shapes (Revised edition)}}, Hemisphere
  Publishing Corp., New York, 1987.
\newline\urlprefix\url{http://adsabs.harvard.edu/cgi-bin/nph-bib\_query?bibcode=1987hpc..book.....B}

\bibitem{GSL}
\href{https://www.gnu.org/software/gsl/}{{GSL --- GNU Scientific Library}}.
\newline\urlprefix\url{https://www.gnu.org/software/gsl/}

\end{thebibliography}
\else
\providecommand{\noopsort}[1]{}\providecommand{\singleletter}[1]{#1}%

\fi

\end{document}